\documentclass{svjour3} 
\usepackage{graphicx} 
\usepackage{hyperref}
\usepackage{amsmath,bm}
\usepackage{microtype}
\usepackage{array}
\usepackage{float}
\usepackage{xurl}
\usepackage{amsmath}
\usepackage{listings}
\usepackage{mathtools}
\usepackage{amsfonts}
\usepackage{algpseudocode}
\usepackage{color}
\usepackage{algorithm}
\RequirePackage{fix-cm}
\smartqed  
\usepackage{graphicx}
\usepackage{cite}

\usepackage{mathptmx}      
\usepackage{latexsym}
\usepackage{subfigure}
\usepackage{amsmath}
\usepackage{booktabs}
\usepackage{float}
\usepackage{multicol}
\usepackage{multirow}
 \usepackage{array}
\usepackage{mathtools}
\usepackage{amsfonts}
\usepackage{amssymb}
\usepackage{algpseudocode}
\usepackage[misc]{ifsym}
\usepackage{color}
\usepackage{algorithm}

\usepackage{graphicx}
\graphicspath{ {./images/} }
\begin{document}

\title{Adaptive multi-gradient methods for quasiconvex vector optimization and applications to multi-task learning}
	
	\titlerunning{Adaptive multi-gradient methods for quasiconvex vector optimization}        
	
	\author{Nguyen Anh Minh \and Le Dung Muu \and Tran Ngoc Thang
	}
	
	
	\institute{
		Nguyen Anh Minh \at
		Faculty of Mathematics and Informatics, Hanoi University of Science and Technology, 1st Dai Co Viet street, Hanoi, Viet Nam \\
		\email{minh.na194117@hust.edu.vn}  \\
 		Le Dung Muu \at
		TIMAS, Thang Long University and Institute of Mathematics, VAST, Hanoi, Vietnam \\
		\email{ldmuu@math.ac.vn} \\
  		Tran Ngoc Thang (\;\Letter\;) \at
		Faculty of Mathematics and Informatics, Hanoi University of Science and Technology, 1st Dai Co Viet street, Hanoi, Viet Nam \\
		\email{thang.tranngoc@hust.edu.vn}  \\
	}
	
	\date{\today}

	\maketitle
	
	\begin{abstract}
We present an adaptive step-size method, which does not include line-search techniques, for solving a wide class of nonconvex multiobjective programming problems on an unbounded constraint set. We also prove convergence of a general approach under modest assumptions. More specifically, the convexity criterion might not be satisfied by the objective function. Unlike descent line-search algorithms, it does not require an initial step-size to be determined by a previously determined Lipschitz constant. The process's primary characteristic is its gradual step-size reduction up until a predetermined condition is met.
It can be specifically applied to offer an innovative multi-gradient projection method for unbounded constrained optimization issues. Preliminary findings from a few computational examples confirm the accuracy of the strategy. We apply the proposed technique to some multi-task learning experiments to show its efficacy for large-scale challenges. 

		
	\end{abstract}
\section{Introduction}\label{s1}
	
With a multitude of real-world applications, gradient descent methods are a popular tool for a broad variety of programming problems, including both scalar and vector optimization, from convex to nonconvex, and references \cite{Boyd}, \cite{Cevher}, and \cite{Lan}. Iterative solutions based on gradient directions and step sizes are produced by gradient descent algorithms at each iteration.
The step-size was chosen using one of the few widely used methods, and for a long time, researchers concentrated on determining the direction to enhance the convergence rate of techniques (see \cite{Boyd}, \cite{Nesterov} ).


In order to lower the total computational cost of the method, new key areas of machine learning applications with large dimensionality and nonconvex objective functions have necessitated the invention of novel step-size choosing procedures (see \cite{Cevher}, \cite{Lan}). One-dimensional minimization line-search, whether precise or approximate, is computationally expensive every iteration, especially when determining the function value is almost the same as determining its derivative and necessitates solving intricate auxiliary issues (see \cite{Boyd}). The step-size value can be determined using known information, such as the gradient's Lipschitz constants, to skip the line-search. But doing so necessitates only employing a portion of their imprecise estimations, which slows convergence. For the well-known divergent series rule, this is also true (see \cite{Kiwiel}, \cite{Nesterov}).

There are numerous applications of the class of quasiconvex multiobjective functions in real-world issues, such as those in economics.
An further drawback of the weighting method for these kinds of issues is that positive combinations of quasiconvex functions might not be quasiconvex. The class of quasiconvex problems has been the subject of several research (\cite{Luc2005}, \cite{LucBook}). For the purpose of minimizing a quasiconvex-concave function under convex and quasiconvex-concave inequality constraints, a class of branch-and-bound algorithms is put forth \cite{Horst}. Despite the fact that Kiwiel created the gradient approach for quasiconvex programming in 2001 \cite{Kiwel}, the usage of decreasing step-size causes the approach to converge slowly.  After that, the gradient approach has been improved. Examples of these include the work of Yu et al. (\cite{Yu}, 2019) and Hu et al. (\cite{Hu}, 2020), who both utilize a constant step-size but require that the objective function meet the Holder condition. The other technique is the neurodynamic approach (see \cite{Bian} (Bian et al, 2018), \cite{Liu} (Liu et al., 2022)) that solves pseudoconvex programming problems with unbounded constraint sets using recurrent neural network models. It does not select the step size in an adaptive manner.

The adaptive step-size procedure was proposed in \cite{Ferreira} (Ferreira et al., 2020) and \cite{Konnov} (Konnov, 2018) and is mostly used to solve single-objective optimization problems. The step-size algorithm for the conditional gradient approach is proposed in \cite{Konnov}, and this method works well for solving pseudoconvex programming problems where the feasible set is bounded. The scope of \cite{Ferreira} has been extended to include the specific scenario of an unbounded feasible set, which is not applicable in the unconstrained scenario. 

Adaptive step sizes have been used in some recent research to address multi-objective optimization problems; \cite{ghalavand} and \cite{assuncao} are two noteworthy examples. In \cite{ghalavand}, unconstrained multi-objective optimization problems were solved by combining the BFGS quasi-Newton approach with a nonmonotone adaptive step size strategy. Additionally, a conditional gradient method combined with an adaptive step size strategy was used in \cite{assuncao} to solve limited multi-objective optimization problems. Their adaptive step size, however, is only theoretical because it was created using knowledge of the objective function's Lipschitz constant.

In this study, we present a novel monotone adaptive step-size technique that does not require line searches for a wide range of unconstrained multiobjective optimization problems in which the objective function is smooth and nonconvex. Reducing the step size gradually until a predefined condition is met is an essential part of this process.Additionally, the Lipschitz continuity of the objective function's gradient is not necessary for convergence; the Lipschitz constant is only used in one particular algorithmic scenario to determine the step size. Preliminary computer tests demonstrate the effectiveness of the proposed modification. We demonstrate the effective performance of the proposed method on large-scale problems by conducting a series of machine learning experiments, such as multi-variable logistic regressions and neural networks for classification.

In this work, we present a new adaptive step-size approach that does not require line searches for a wide class of multiobjective programming problems in which the constraint set is unbounded closed convex and the objective function is nonconvex smooth. Reducing the step size gradually until a predefined condition is met is an essential part of this process. The method does not use predefined constants, even though convergence requires the Lipschitz continuity of the objective function's gradient. Preliminary computer experiments have demonstrated the effectiveness of the proposed modification. We conduct a range of machine learning experiments, such as multi-task learning and classification using neural networks, to demonstrate the effectiveness of the suggested approach on large-scale problems.

Chapter 3 will present our Adaptive Multi-gradient algorithm for solving Multi-objective optimization Problems while proving the convergence of the algorithm. Chapter 4 will provide an upgraded version of the Adaptive Multi-gradient algorithm to help evenly distribute the found optimal solutions on the Pareto Front. Chapter 5 will provide some examples of solving multi-objective optimization problems for convex, pseudoconvex, and nonconvex functions. Chapter 6 will discuss the application of the proposed algorithm in solving Multi-task Learning problems with specific tasks in image processing and text processing. Chapter 7 will present some conclusions of this paper.
	\section{Preliminaries and Notations}\label{s2}

We examine a scenario involving a single-objective function $f:\mathbb{R}^m \to \mathbb R$ with a permissible set $U\subseteq \mathbb{R}^m$. The subsequent definitions are derived from \cite{Mangasarian}.

\begin{definition}\label{def1}\cite{Mangasarian}
The differentiable function $f:\mathbb{R}^m \to \mathbb R$ is considered
\begin{itemize}
\item[i)] convex on $ U $ if, for all $ x, y \in U $, and $ \lambda \in [0,1] $, the following inequality holds:
\begin{equation*}
f(\lambda x +(1-\lambda)y)\leq \lambda f(x)+(1-\lambda)f(y).
\end{equation*}
\item[ii)] pseudoconvex on $ U $ if, for all $ x, y \in U $, the condition is met:
\begin{equation*}
\langle \nabla f(x), y-x \rangle \geq 0 \Rightarrow f(y)\geq f(x).
\end{equation*}
\item[iii)] quasiconvex on $ U $ if, for all $ x, y \in U $, and $ \lambda \in [0;1] $, the inequality holds:
\begin{equation*}
f(\lambda x +(1-\lambda)y)\leq \max \{f(x);f(y) \}.
\end{equation*}
\end{itemize}
\end{definition}
	\begin{proposition}\label{P_1}\cite{Dennis}
The differentiable function $ f $ is quasiconvex on $ U $ if and only if
		\[ 
		f(y)\leq f(x) \Rightarrow 		\langle \nabla f(x), y-x \rangle \leq 0.
		\] 
	\end{proposition}
	It is important to note that the implication chain holds: "$f$ is convex" $\Rightarrow$ "$f$ is pseudoconvex" $\Rightarrow$ "$f$ is quasiconvex" \cite{Mangasarian}.

The range, or image space, of a matrix $M \in \mathbb{R}^{m \times n}$ is denoted by $R(M)$, and $I \in \mathbb{R}^{n \times n}$ represents the unit matrix. For two matrices $A, B \in \mathbb{R}^{n \times n}$, $B \leq A$ ($B < A$) signifies that $A-B$ is positive semidefinite (definite). In the upcoming discussions, the Euclidean norm in $\mathbb{R}^n$ is denoted by $\|\cdot\|$, and $B[x, r]$ represents the closed ball of radius $r$ with center $x \in \mathbb{R}^n$. The same notation $\|\cdot\|$ is employed for the induced operator norms on the corresponding matrix spaces.

Define, for $A \in \mathbb{R}^{m \times n}$
$$
\|A\|_{\infty, 2}:=\max _{x \neq 0} \frac{\|A x\|_{\infty}}{\|x\|} .
$$
Then $\|\cdot\|_{\infty, 2}$ is a norm in $\mathbb{R}^{m \times n}$. It is easy to prove that
$$
\begin{aligned}
\|A\|_{\infty, 2} & :=\max _{i=1, \ldots, m}\left\|A_{i, \cdot}\right\| \\
& =\max _{i=1, \ldots, m}\left(\sum_{j=1}^n A_{i, j}^2\right)^{1 / 2} .
\end{aligned}
$$

\subsection{Multi-objective optimization problem}\label{F_Lipschitz}
We examine the multi-objective optimization problem represented as follows:
\begin{equation}\label{MOP}
\mathrm{Min}_{x\in U} F(x) \tag{MOP($ F,U $)},
\end{equation}
where the set $U \subseteq \mathbb{R}^n$ is nonempty, convex, and closed. The vector function $F: U \longrightarrow \mathbb{R}^m$ is assumed to be differentiable on an open set containing $U$. We will assume that $U = \{x\in\mathbb{R}^n|g_i(x)\leq 0, i=1,\ldots,l \}$. The Jacobian matrix of $F$ at $x$ is denoted as follows:
$$
JF(x)=\left[\begin{array}{c}
\nabla F_1(x)^T \\
\vdots \\
\nabla F_m(x)^T
\end{array}\right] \in \mathbb{R}^{m \times n}
$$
In this paper, we consider the Jacobian matrix $JF$ of $F$ to be continuously Lipschitz in $U$. This implies that the component-wise derivatives $\nabla F_j \in \mathbb{R}^n$ of the function $F_j$ are continuously Lipschitz for all $j=1,\ldots,m$. It is also assumed that the solution set of \eqref{MOP} is nonempty. The subsequent definitions are sourced from \cite{ehrgott2005multicriteria}.

\begin{definition}  An element $x^* \in U$ is referred to as a \textit{efficient point} or a Pareto optimal point of problem (\ref{MOP}) if there is no $y \in U, F(y) \leq F\left(x^*\right)$ and $F(y) \neq F\left(x^*\right)$, where the $\leq$ symbol between vectors is interpreted as the comparison between each component function. There is no $y \in U$ such that $F_i(y) \leq F_i\left(x^*\right)$ for every $i=1, \ldots, m$, and there exists $i_0 \in{1, \ldots, m}$ such that $F_{i_0}(y)<F_{i_0}\left(x^*\right)$.The Pareto Front is defined as the collection of Pareto optimum locations.
\end{definition}

Given the set $\mathbb{R}_{+}^m :=\mathbb{R}_{+} \times \cdots \times \mathbb{R}_{+}$, the partial order between vector function values is defined as follows: $$ F(y) \leq F(z) \Longleftrightarrow F(z)-F(y) \in \mathbb{R}_{+}^m, $$ Hence, in order to solve problem (\ref{MOP}), it is necessary to find a point that represents the minimum in the partial order.

\begin{definition} A point $x^* \in U$ is referred to as a \textit{weakly efficient} solution or a \textit{weak Pareto optimal point} of problem (\ref{MOP}). If there is no element in U such that F(y) is less than F(x*), where the inequality F(y) < F(x*) is interpreted as a partial order.
The definition ends.
\end{definition}

Considering $\mathbb{R}_{++}^m :=\mathbb{R}_{++} \times \cdots \times \mathbb{R}_{++}$, we also define

\begin{align*}
-\mathbb{R}_{++}^m=\{-s: s \in \mathbb{R}_{++}^m\}.
\end{align*}

\begin{definition} A point $x^* \in U$ is considered locally efficient, or weakly locally efficient, in problem (\ref{MOP}) if there exists a neighborhood $V \subseteq U$ of $x^*$ where $x^*$ is an efficient, or weakly efficient, point in problem (\ref{MOP}) when restricted to $V$.
\end{definition}

If the set $U$ is convex and the function $F$ is $\mathbb{R}_{+}^m$-convex (meaning that each component function of $F$ is convex), then every locally efficient point is also a globally efficient point.
Furthermore, every point that is locally efficient is also a weakly locally efficient point.

\begin{definition}\cite{fliege2000steepest}
A vector $z$ belonging to the set $U$ is called a tangent direction of $U$ at $z$ if there exists a sequence $\left(z^{k}\right)_{k}$ contained in $U$ and a positive scalar $\lambda \in \mathbb{R}^{+}$ such that $$ \lim _{k \rightarrow \infty} z^{k}=z \quad \text { and } \quad \lim _{k \rightarrow \infty} \lambda \frac{z^{k}-z}{\left\|z^{k}-z\right\|}=s $$
The collection of all tangent directions of the set $U$ at the point $z$ is referred to as the tangent cone of $U$ at $z$ and symbolized as $T(U, z)$.

\end{definition}
For a closed convex set $U$, we have $T(U, z) := U(z) := \{s \in \mathbb{R}^n | s = u - z \text{ for some } u \in U \}$.

\begin{definition} \label{Pareto_sta}
A point $x^*$ in the set $U$ is referred to as a \textit{Pareto stationary point} (or \textit{Pareto critical point}) of the vector function $F$ over $U$ if the intersection of the product of the Jacobian matrix $JF(x^*)$ and the set $U(x^*)$ with the negative real numbers raised to the power of $m$ is empty.
\end{definition}
This definition is a necessary condition (but not sufficient) for a point to be an effective Pareto point and was first used in \cite{fliege2000steepest} to define the steepest descent algorithm for multiobjective optimization.

 \begin{lemma} \label{F_pseudo} If the function $F$ is pseudo-convex, meaning that its component functions $F_i, i=1,\ldots,m$ are pseudo-convex, then Definition \ref{Pareto_sta} serves as both a necessary and sufficient condition for a point to be considered a weakly efficient point.
\end{lemma}

Suppose, on the contrary, that there is a point $y \in U$ such that $F(y) < F(x)$. Given that $F_i$ is pseudo-convex for each $i=1, \ldots, m$, we may conclude that $\left\langle\nabla F_i(x),(y-x)\right\rangle<0$. Consequently, $JF(x)(y-x) \in-\mathbb{R}_{++}^m$, which violates Definition \ref{Pareto_sta} (it is important to note that we are assuming the set $U$ to be both convex and closed).

\subsection{Multi-task Learning Problem}\label{MTL_MOP}
Traditional machine learning methods typically focus on solving one task with a single model. This approach may cause us to overlook valuable information that could help us perform better on the task at hand, information that can be derived from related tasks. To illustrate this, consider a simple example where you want to predict the price of a house based on features such as area, number of rooms, number of floors, proximity to commercial centers, and so on. Clearly, adding a few additional tasks, such as predicting whether the house is in the city or suburbs or whether it's a villa or an apartment, can provide valuable insights for predicting its price.

From a specialized perspective, Multi-Task Learning (MTL) is often implemented by sharing common model parameters (shared parameters) among network architectures and specific parameters for each individual task's network architecture. This allows us to obtain a more generalized model for the original task we are interested in.

Our research effort aims to optimize shared parameters using multi-objective optimization methods. Multi-Task Learning has been successfully applied in computer vision (Bilen and Vedaldi, 2016; Misra et al., 2016; Kokkinos, 2017; Zamir et al., 2018), natural language processing (Collobert and Weston, 2008; Dong et al., 2015; Liu et al., 2015a; Luong et al., 2015; Hashimoto et al., 2017), and audio processing (Huang et al., 2013; Seltzer and Droppo, 2013; Huang et al., 2015). It is also referred to by other names such as Joint Learning, Learning to Learn, and Learning with Auxiliary Tasks. When we encounter an optimization problem that involves many loss functions, we are really dealing with a problem that is connected to Multi-Task Learning.

The MTL problem is represented as a Multi-Objective Programming (MOP) problem, using the approach proposed by Sener et al. \cite{sener}. The MTL problem entails the execution of $m$ jobs using a loss function represented by a vector.
The equation \ref{MTL_MOP} represents a multi-task learning problem, where we aim to minimize the loss function $\mathcal{L}(\theta)$ with respect to the parameter set $\theta$. The loss function $\mathcal{L}_i(\theta)$ corresponds to the loss function for task $i$ over all tasks. MTL algorithms concurrently optimize many tasks. Problem (\ref{MTL_MOP}) is a multi-objective optimization problem, and often, there is no single solution or parameter configuration that can maximize all tasks concurrently. Thus, we can exclusively offer Pareto optimal solutions that illustrate the compromise between several tasks.

\section{Adaptive multi-gradient for solving Multi-objective optimization problem}
\subsection{Solving subproblem by KKT condition}
We will now solve issue \eqref{MOP} when $U$ is the set of all real $n$-dimensional vectors. The steepest descent direction $s(x)$ at $x$ is defined as the answer to the optimization problem:
\begin{equation}\label{3.5}
    \left\{\begin{array}{cl}
\min & \max _{j=1, \ldots, m} \nabla F_j(x)^T s+\frac{1}{2} ||s||^2 \\
\text { subject to } & s \in \mathbb{R}^n,
\end{array}\right.
\end{equation}

The problem \eqref{3.5} possesses a solitary optimal solution due to the strong convexity of the function $\nabla F_j(x)^T s+\frac{1}{2} s^Ts$ with respect to the variable $s$ for $j=1, \ldots, m$. When $m=1$, the steepest descent path $s(x)$ is equal to the negative gradient of $F(x)$, denoted as $-\nabla F(x)$.

The formulation of problem \eqref{3.5} is based on the idea of approximating the expression $$ \max _{j=1, \ldots, m} F_j(x+s)-F_j(x) $$ by employing a first-order Taylor expansion around $x$ for each function $F_j$. The authors in \cite{fliege2000steepest} have suggested a formulation of the optimization problem \eqref{3.5} to determine the steepest descent direction for multi-objective optimization. The most favorable solution for equation \eqref{3.5} will be represented as $\Theta(x)$. Thus, we have the following equations: 
\begin{equation}\label{3.6}
    \Theta(x)=\inf _{s \in \mathbb{R}^n} \max _{j=1, \ldots, m} \nabla F_j(x)^T s+\frac{1}{2} ||s||^2,
\end{equation}
and
\begin{equation} \label{3.7}
    s(x)=\arg \min _{s \in \mathbb{R}^n} \max _{j=1, \ldots, m} \nabla F_j(x)^T s+\frac{1}{2} ||s||^2
\end{equation}
Despite the non-smooth nature of issue \eqref{3.5}, it may still be classified as a second-order optimization problem and effectively solved using the Karush-Kuhn-Tucker (KKT) criteria. The problem \eqref{3.5} may be expressed as the following optimization problem, denoted as ($P_1$):
\begin{equation}\label{projected}\tag{$P_1$}
\begin{cases}\min & h(t, s)=t \\ \text { subject to } & \nabla F_j(x)^T s+\frac{1}{2} ||s||^2-t \leq 0 \quad(1 \leq j \leq m) \\ & (t, s) \in \mathbb{R} \times \mathbb{R}^n.\end{cases}
\end{equation}
So the Lagrangian objective function is
$$
L((t, s), \lambda)=t+\sum_{j=1}^m \lambda_j\left(\nabla F_j(x)^T s+\frac{1}{2} ||s||^2-t\right) .
$$
The problem denoted by equation \eqref{projected} possesses a solitary optimal solution, denoted as $(\Theta(x), s(x))$, due to its convex nature and the fulfillment of Slater's condition. According to the Karush-Kuhn-Tucker (KKT) criteria, there is a KKT multiplier vector $\lambda=\lambda(x)$ that corresponds to $s=s(x)$ and $t=\Theta(x)$. We can readily calculate the expression: \begin{equation}\label{3.13_no} s(x)=-\sum_{j=1}^m \lambda_j(x) \nabla F_j(x) . \end{equation}

The most steep descent direction is determined by the weight vector $\lambda$ acquired from solving the KKT system, which converts the multi-objective optimization issue into a single-objective optimization problem. If the regularity criterion is met, the KKT solution for the convex programming problem \eqref{projected} corresponds to the solution of both the primal and dual problems. When we combine it with equation \eqref{3.13_no}, we obtain:

\begin{equation}\label{3.14_no}
\begin{split}
     \Theta(x) &=\sup _{\lambda \geq 0} \inf _{s \in \mathbb{R}^n} L((t, s), \lambda) \\
     &=\sup _{\substack{\lambda \geq 0 \\ \sum \lambda_j=1}} \inf _{s \in \mathbb{R}^n} \sum_{j=1}^m \lambda_j\left(\nabla F_j(x)^T s+\frac{1}{2} s^Ts\right) \\
     &=\sup _{\substack{\lambda \geq 0 \\ \sum \lambda_j=1}} \inf _{s \in \mathbb{R}^n} \left(\sum_{j=1}^m \lambda_j(x)\nabla F_j(x)^T. s+\frac{1}{2} s^Ts \right)\\
     & = \sup _{\substack{\lambda \geq 0 \\ \sum \lambda_j=1}}  -\dfrac{1}{2}\|\sum_{j=1}^m\lambda_j(x)\nabla F_j(x)\|^2
\end{split}
\end{equation}

In the following part, we will analyze the attributes of $\Theta(x)$ and investigate its relationship with $s(x)$. Moreover, we examine the stationarity of $\Theta(x)$. The subsequent lemmas and remark are results reported in the publication referenced as \cite{fliege2000steepest}.

\begin{lemma}\label{lemma4.2}
    With $x \in \mathbb{R}^n$, we have
$$
\Theta(x)=-\frac{1}{2} s(x)^Ts(x).
$$
\end{lemma}

\begin{lemma}\label{lemma 3.2_no} 
Denote $s(x)$ as the solution and $\Theta(x)$ as the optimum value of problem \eqref{projected}.
\begin{enumerate}
    \item If $x$ is Pareto critical, then $s(x)=0 \in \mathbb{R}^n$ and $\Theta(x)=0$.
    \item If $x$ is not Pareto critical, then $\Theta(x)<0$.
\item The mappings $x \mapsto s(x)$ and $x \mapsto \Theta(x)$ are continuous.

\end{enumerate}
\end{lemma}

Instead of precisely solving problem \eqref{projected}, it is intriguing from an algorithmic standpoint to handle approximate answers. If $x$ is not Pareto critical, we define $s$ as an approximate solution to equation \eqref{projected} with a tolerance $\sigma \in] 0,1]$ if
\begin{equation}\label{Theta<0}
   \max _j(JF(x) s)_j+\frac{1}{2}\|s\|^2 \leq \sigma \Theta(x) \leq 0, \quad s \in \mathbb{R}^n 
\end{equation}
where $\Theta(x)$ is the optimum value of problem \eqref{projected}. 

\begin{lemma}\label{bichan}
Let's assume that $x$ is not Pareto critical and that $s$ is an estimated solution of equation \eqref{projected} with a tolerance value $\sigma \in] 0,1]$. Then
$$
\|s(x)\| \leq 2\|JF(x)\|_{\infty, 2}
$$   
\end{lemma}

\textbf{Proof. } If $x$ is not Pareto-critical $\Rightarrow \Theta(x)<0$. Also, $\Theta$ is an approximate solution of \eqref{projected} so
$$
\begin{aligned}
& f_x(s)+\frac{1}{2}\|s\|^2 \leq \sigma \Theta(x) \\
 \Rightarrow\|s\|^2 & \leq 2(\sigma\Theta(x)-f_x(s)) \leq 2(-f_x(s)) \\
& \leqslant 2(-J F(x) s)_j \quad \forall j=\overline{1, m} \\
& \leq 2\|J F(x) s\|_{\infty, 2} \\
& \leq 2 \| J F\left(x)\left\|_{\infty, 2}\right\| s \|\right. \\
\end{aligned}$$

\subsection{Compute the adaptive step length}
Assume we have a vector $s \in \mathbb{R}^n$ with $J F(x) s<0$. To compute a step length $\alpha>0$ we use an adaptive rule.  Let $\sigma$ be a constant that is chosen in the interval $(0,1)$. The criterion for accepting $\alpha$ is \begin{equation}\label{steplength}
    F(x+\alpha s) \leq F(x)+\sigma \alpha J F(x) s .
\end{equation}

Let us begin with an initial value of $\alpha$ that lies between 0 and 1. If the condition is not met, we update $\alpha$ by multiplying it with $\kappa$, where $\kappa$ is a value between 0 and 1 (inclusive). We then go to the next iteration and obtain a new direction, denoted as $s$. The finiteness of this technique can be deduced from the fact that equation \eqref{steplength} is strictly valid for sufficiently tiny values of $\alpha>0$.

\begin{lemma}\cite{fliege2000steepest}\label{JF(x)v<0}
  If $F$ is differentiable and $J F(x) s<0$, then there exists a positive value $\varepsilon$ (which may depend on $x, s$, and $\sigma$) such that the inequality $$ F(x+t s)<F(x)+\sigma \alpha J F(x) s $$ holds for any $\alpha$ in the interval $(0, \varepsilon]$.
\end{lemma}


\subsection{The complete algorithm}
At each step, for a non-stationary point, we solve the problem \eqref{projected} to determine the steepest descent direction using the formula in Lemma \ref{lemma4.2}. We proposed the following algorithm to solve problem \eqref{MOP} with $U=\mathbb{R}^n$.

\begin{algorithm}[H]
\caption{Adaptive multi-gradient}
\label{algo_nocons_1}
\textbf{Step 1}. Set $\kappa \in (0,1]$, $\sigma \in [0,1]$, $\alpha_1 \in (0, 1)$, and $x^{1} \in \mathbb{R}^n$. Initialize $k:=1$.


\textbf{Step 2.} (Main loop)

(a) Find the descent direction by solving the following problem:

$$
\begin{aligned}
& \min _{\lambda_j}    ||\sum_{j=1}^m \lambda_j^k\nabla F_j(x^k)||^2 \\
&\text {subject to } \lambda_j \geq 0,\,\sum_{j=1}^m \lambda_j^k=1
\end{aligned}
$$

From this, we obtain the optimal solution:

$\mathbf{\lambda^k} = (\lambda_j^k)$,$\,\,s(x^k)=-\left(\sum_{j=1}^m \lambda^{k}_{j} \nabla F_j(x^{k})\right)$ and $\Theta(x^{k}) =-\dfrac{1}{2}||s(x^k)||^2.$ 

(b) If $\Theta(x^{k})=0$, then \textbf{End}. Otherwise, continue to \textbf{Step 2(c).}

(c) Set $x^{k+1}:=x^{k}+\alpha_{k} s(x^{k})$.

(d) Compute the step size:

If \begin{equation}\label{kiemtra}
F(x^{k+1}) \leq F(x^{k}) +\sigma \left\langle JF(x^k),x^{k+1}-x^k \right\rangle
\end{equation}

then $\alpha_{k+1} = \alpha_k$, otherwise set $\alpha_{k+1}:=\kappa \alpha_k$.

\textbf{Step 3.} Set $k:=k+1$ and go to \textbf{Step 2.}
\end{algorithm}
\subsection{Convergence of adaptive multi-gradient methods}
The convergence analysis of the proposed algorithm is based on quasi-Fejer convergence. We recall that a sequence $\left\{z^k\right\} \subset \mathbb{R}^n$ is said to be quasi-Fejér convergent to a set $U, U \neq \emptyset$, if and only if for each $z \in U$, there exists a sequence $\left\{\epsilon_k\right\} \subset \mathbb{R}_{+}$such that $\sum_{k=1}^{+\infty} \epsilon_k<+\infty$ and
$$
\left\|z^{k+1}-z\right\|^2 \leq\left\|z^k-z\right\|^2+\epsilon_k \text {. }
$$
\begin{theorem}\label{theorem1_nocons}
		Assume that the sequence $ \left\{ x^k\right\} $ is generated by Algorithm \ref{algo_nocons_1}. Then, each limit point (if any) of the sequence $ \left\{ x^k\right\} $ is a Pareto stationary point of the problem \eqref{MOP} where $U=\mathbb{R}^n$. 
	\end{theorem}
\textbf{Proof.} Let $y$ be an accumulation point of the sequence $\left(x^k\right)$ and let $s(y)$ and $\Theta(y)$ be the solution and the optimum value of \eqref{projected} at $y$.

According to Lemma \ref{lemma 3.2_no} it is enough to prove that $\Theta(y)=0$. Clearly, the sequence $F\left(x^k\right)$ is componentwise strictly decreasing and we have
$$
\lim _{k \rightarrow \infty} F\left(x^k\right)=F(y)
$$
Therefore,
$$
\lim _{k \rightarrow \infty}\left\|F\left(x^k\right)-F\left(x^{k+1}\right)\right\|=0 .
$$

Let's consider the sequence $\{k_p\} = \{k_1, k_2, ...\}$ with $p \rightarrow \infty$, satisfying the inequality \eqref{kiemtra}. In that case, the value of $\alpha_{k_p}$ will remain unchanged.

\textbf{Case 1:} $\{k_p\}$ is finite.

If the sequence $\{k_p\}$ is finite with sufficiently large $p$, there exists a step $K$ onwards where the inequality \eqref{kiemtra} is not satisfied. In other words,
$$
F_j\left(x^{K+1}\right) \geq F_j(x^{K})+\sigma \alpha_{K}
\langle J F(x^{K}), s(x^{K})\rangle_j
$$
for at least one $j \in\{1, \ldots, m\}$.
According to Proposition \ref{JF(x)v<0}, we have $\max_j\langle J F(x^{K}), s(x^{K})\rangle \geq 0$. Combining this with \eqref{Theta<0}, we obtain $\Theta(x^{K})=0$, which implies that $x^K$ is a Pareto stationary point.

\textbf{Case 2:} $\{k_p\}$ is infinite.

In the infinite sequence $\{k_p\}$, we will always have
$$
F\left(x^k\right)-F\left(x^{k+1}\right) \geq -\alpha_k \sigma \langle JF(x^k), s(x^k) \rangle \geq 0,
$$
for all $k\in \{k_p\}$ and therefore
\begin{equation}\label{gioihan}
\lim _{k \rightarrow \infty} \alpha_k \sigma \langle JF(x^k), s(x^k) \rangle=0 .
\end{equation}
Observe that $\left.\left.\alpha_k \in\right] 0,1\right]$ for all $k\in \{k_p\}$. Now take a subsequence $\left(x^{k_u}\right)_u$ converging to $y$. We will consider two possibilities:
$$
\limsup _{u \rightarrow \infty} \alpha_{k_u}>0
$$
and
$$
\limsup _{u \rightarrow \infty} \alpha_{k_u}=0 .
$$

\textbf{Case 2.1} In this case, there exist a subsequence $\left(x^{k_{\ell}}\right)_{\ell}$ converging to $y$ and satisfying
$$
\lim _{\ell \rightarrow \infty} \alpha_{k_{\ell}}=\bar{\alpha}>0
$$
Using \eqref{gioihan} we conclude that
$$
\lim _{\ell \rightarrow \infty} J F\left(x^{k_{\ell}}\right) s^{k_{\ell}}=0
$$
which also implies
$$
\lim _{\ell \rightarrow \infty} \Theta\left(x^{k_{\ell}}\right)=0
$$
Since $x \mapsto \Theta(x)$ is continuous, we conclude that $\Theta(y)=0$, so $y$ is Pareto critical.

\textbf{Case 2.2} Due to the results in Lemma \ref{bichan} we know that the sequence $\left(s^{k_u}\right)_u$ is bounded. Therefore we can take a subsequence $\left(x^{k_r}\right)_r$ of $\left(x^{k_u}\right)_u$ such that the sequence $\left(s^{k_r}\right)_r$ also converges to some $\bar{s}$. Note that for all $r$ we have
$$
\max _i\left(J F\left(x^{k_r}\right) s^{k_r}\right)_i \leq \sigma \Theta\left(x^{k_r}\right)<0
$$
Therefore, passing onto the limit $r \rightarrow \infty$ we get
\begin{equation}\label{Theta}
    \frac{1}{\sigma} \max _i(J F(y) \bar{s})_i \leq \Theta(y) \leq 0
\end{equation}
Take some $q \in \mathbb{N}$. For $r$ large enough,
$$
\alpha_{k_r}<\alpha_1^q
$$
which means that the Adative step size \eqref{kiemtra} is not satisfied for $\alpha=\alpha_1^q$, i.e.,
$$
J F\left(x^{k_r}+\left(\alpha_1^q\right) s^{k_r}\right) \nleq F\left(x^{k_r}\right)+\sigma\left(\alpha_1^q\right) J F\left(x^{k_r}\right) s^{k_r}
$$
(for $r$ large enough). Passing onto the limit $r \rightarrow \infty$ (along a suitable subsequence, if necessary) we get
$$
F_j\left(y+\left(\alpha_1^q\right) \bar{s}\right) \geq F_j(y)+\sigma\left(\alpha_1^q\right)(J F(y) \bar{s})_j
$$
for at least one $j \in\{1, \ldots, m\}$. Note that this inequality holds for any $q \in \mathbb{N}$. From Lemma \ref{JF(x)v<0} it follows that
$$
\max _i(J F(y) \bar{s})_i \geq 0
$$
which, together with \eqref{Theta}, implies $\Theta(y)=0$. So, again we conclude that $y$ is Pareto critical.

Now, we are interested in the study of the convergence properties of the proposed algorithm when the objective function is quasiconvex. Define
$$
T=\left\{x \in \mathbb{R}^n: F(x) \leq F\left(x^k\right), \forall k\right\} .
$$
\begin{proposition}\label{quasi}\cite{cruz2011convergence}
If $F$ is quasiconvex and $x \in \mathbb{R}^n $ then, $\left\langle\nabla F_i\left(x^k\right), x-x^k\right\rangle \leq 0$ for all $k$ and all $i=1, \ldots, m$.
\end{proposition}
\begin{theorem}\label{theorem2}
 Assume that $F$ is a quasiconvex function on $\mathbb{R}^n$. If $T \neq \emptyset$ then $\{x^k\}$ is generated by Algorithm \ref{algo_nocons_1} converges to a stationary point. 
\end{theorem}
\textbf{Proof.} Based on the proof of Theorem \ref{theorem1_nocons}, we again consider the sequence $\{k_p\} = \{k_1, k_2, ...\}$ with $p \rightarrow \infty$, satisfying the inequality \eqref{kiemtra}. 

\textbf{Case 1:} $\{k_p\}$ is finite.

Due to the finiteness of the sequence $\{k_p\}$, we can prove the existence of a Pareto critical point $x^K$ after a finite number of steps. Therefore when $F$ is a quasiconvex function on $\mathbb{R}^n$, $\{x^k\}$ converges to a stationary point.

\textbf{Case 2:} $\{k_p\}$ is infinite.

After some simple algebra, for each $x \in \mathbb{R}^n$, we obtain
$$
d_k:=\left\|x^k-x\right\|^2-\left\|x^{k+1}-x\right\|^2+\left\|x^{k+1}-x^k\right\|^2=2\left\langle x^{k+1}-x^k, x-x^k\right\rangle .
$$

Using the update formula $x^{k+1}=x^k+\alpha_k s^k$, it follows that
$$
\begin{aligned}
d_k & =2 \alpha_k\left\langle s^k, x-x^k\right\rangle \\
& =2 \alpha_k\left\langle s^k, x-x^k-s^k+s^k\right\rangle \\
& =2 \alpha_k\left\langle s^k, x-x^k-s^k\right\rangle+2 \alpha_k\left\|s^k\right\|^2 .
\end{aligned}
$$
So we get
$$
d_k-2 \alpha_k\left\|s^k\right\|^2  =2 \alpha_k\left\langle -\sum_{j=1}^m \lambda_j^k \nabla F_j\left(x^k\right), x-x^k-s^k\right\rangle
$$
Therefore,
$$
\begin{aligned}
\left\|x^k-x\right\|^2-\left\|x^{k+1}-x\right\|^2+\alpha_k\left(\alpha_k-2\right)\left\|s^k\right\|^2 & =d_k-2 \alpha_k\left\|s^k\right\|^2 \\
& = 2 \alpha_k \left(\left\langle\sum_{j=1}^m \lambda_j^k \nabla F_j\left(x^k\right), x^k-x\right\rangle+\left\langle\sum_{j=1}^m \lambda_j^k \nabla F_j\left(x^k\right), s^k\right\rangle\right) .
\end{aligned}
$$
Rewriting this last equality and using that $F\left(x^{k+1}\right) \leq F\left(x^k\right)+\sigma \alpha_k JF\left(x^k\right) s^k$ with $\alpha_k<1$, it follows that
$$
\left\|x^{k+1}-x\right\|^2 \leq\left\|x^k-x\right\|^2+2 \alpha_k \left\langle\sum_{j=1}^m \lambda_j^k \nabla F_j\left(x^k\right), x-x^k\right\rangle+\frac{2}{\sigma} \sum_{j=1}^m \lambda_j^k\left(F_j\left(x^k\right)-F_j\left(x^{k+1}\right)\right) .
$$
Using that $F\left(x^{k+1}\right) \leq F\left(x^k\right)$ and $1 \geq \lambda_j^k \geq 0$ for all $j=1, \ldots, m$, we obtain
$$
\sum_{j=1}^m\left(F_j\left(x^k\right)-F_j\left(x^{k+1}\right)\right) \geq \sum_{j=1}^m \lambda_j^k\left(F_j\left(x^k\right)-F_j\left(x^{k+1}\right)\right)
$$
So we have for all $ x \in \mathbb{R}^n \text { and each } k \text {, there exists }\left\{\lambda_j^k\right\}_{j=1}^m \subset[0,1] \text { satisfying } \sum_{j=1}^m \lambda_j^k=1 \text {, and }$
\begin{equation}\label{bđt}
  \left\|x^{k+1}-x\right\|^2 \leq\left\|x^k-x\right\|^2+2 \alpha_k \left\langle\sum_{j=1}^m \lambda_j^k \nabla F_j\left(x^k\right), x-x^k\right\rangle+\frac{2}{\sigma} \sum_{j=1}^m\left(F_j\left(x^k\right)-F_j\left(x^{k+1}\right)\right) .  
\end{equation}
for $k \in \{k_p\}$. We have the sequence $\left\{F\left(x^k\right)\right\}$ is strictly monotone decreasing. Therefore $\varepsilon_k:=\sum_{j=1}^m\left[F_j\left(x^k\right)-F_j\left(x^{k+1}\right)\right]$ is positive for all $k$. 
It follows from Proposition \ref{quasi} that
$$
\left\|x^{k+1}-x\right\|^2 \leq\left\|x^k-x\right\|^2+\frac{2}{\sigma} \varepsilon_k
$$
Observe that the series $\sum_{k=1}^{\infty} \varepsilon_k$ is convergent. Indeed, $\varepsilon_k>0$ for all $k$ and
$$
\begin{aligned}
\sum_{k=0}^{\ell} \varepsilon_k & =\sum_{k=0}^{\ell}\left(\sum_{i=1}^m\left(F_i\left(x^k\right)-F_i\left(x^{k+1}\right)\right)\right)=\sum_{i=1}^m\left(\sum_{k=0}^{\ell}\left(F_i\left(x^k\right)-F_i\left(x^{k+1}\right)\right)\right) \\
& =\sum_{i=1}^m\left(F_i\left(x^0\right)-F_i\left(x^{\ell+1}\right)\right) \leq \sum_{i=1}^m\left(F_i\left(x^0\right)-F_i(x)\right) .
\end{aligned}
$$
Therefore, $\left\{x^k\right\}$ is quasi-Fejér convergent to $T$. Using result on quasi-Fejér convergent Lemma 1 in \cite{cruz2011convergence}, the sequence $\left\{x^k\right\}$ is bounded. Let $\bar{x}$ be an accumulation point of $\left\{x^k\right\}$. The sequence $\left\{F\left(x^k\right)\right\}$ is monotone decreasing, which implies that $\bar{x} \in T$. Using Lemma 1 in \cite{cruz2011convergence} we have $\left\{x^k\right\}$ is convergent to $\bar{x}$, which is stationary by Theorem \ref{theorem1_nocons}.

\begin{theorem}\label{theorem3}
Assume that $F$ is a pseudoconvex function on $\mathbb{R}^n$. If $T \neq \emptyset$, then $\{x^k\}$ generated by Algorithm \ref{algo_nocons_1} converges to a weakly efficient point.
\end{theorem}
\textbf{Proof} By \cite{Mangasarian}, pseudoconvexity implies quasiconvexity. Then, using Theorem 2, we
obtain that $\{x_k\}$ converges to a stationary point, which is a weakly efficient point.
\begin{remark}
By using Remark 2 in \cite{fliege2000steepest}, we have the \eqref{kiemtra} criterion holds for $0\leq \lambda \leq 2(1-\sigma)/L$. As a result, if the value of the Lipschitz constant $L$ is already known, we can choose the constant step-size $\lambda\in(0,2/L)$ to solve problem \ref{MOP}.
\end{remark}
\section{Adaptive multi-gradient with preference vectors}
Solving multi-objective problems using scalarization approaches tends to concentrate the obtained Pareto optimal points in a specific region of the Pareto front. In order to evenly distribute the Pareto optimal points across the entire Pareto front, we need to apply additional linear constraints to help partition the image space. On each partitioned image space, we will find a corresponding Pareto optimal solution.

This paper \cite{PMTL} has a similar idea. It looks at a group of $K$ preference vectors $\{\boldsymbol{u}_1, \boldsymbol{u}_2, \ldots, \boldsymbol{u}_K\} \in \mathbb{R}_{+}^m$. We divide the image space into subspaces $\Omega_k (k=1,\ldots,K)$ using the conditions below:
$$
\Omega_k=\left\{\boldsymbol{v} \in R_{+}^m \mid \langle\boldsymbol{u}_i, \boldsymbol{v}\rangle \leq \langle\boldsymbol{u}_k, \boldsymbol{v}\rangle, \forall i=1, \ldots, K\right\}.
$$
The main idea of preference vectors is to decompose a multi-objective problem into
several constrained multi-objective subproblems with
different trade-off preferences among the tasks in
the original problem. By solving these subproblems in
In parallel, the proposed algorithm can obtain a set of well-representative solutions with different trade-offs. Next, we solve the constrained optimization problem on each region $\Omega_k$ by combining the inequality constraints $\mathcal{G}_p$ to aid in the space partitioning.
We need to solve the following problem:
\begin{equation}\label{P3}
    \begin{aligned}
 \min _{x \in \mathbb{R}^n} F(x) & =\left(F_1(x), F_2(x), \cdots, F_m(x)\right) \\
 \text {subject to } & \mathcal{G}_p\left(x\right)=\langle\boldsymbol{u}_p-\boldsymbol{u}_k,F\left(x\right)\rangle \leq 0, \forall p=1, \ldots, K \text {, } 
\end{aligned} 
\end{equation}
It is necessary to sample the preference vectors $u_1, \ldots, u_K$ to cover the objective space evenly to make HV and cosine similarity interact well. This is quite simple in 2D space using $K$ vectors $\left\{\cos \left(\frac{k \pi}{2 K}\right), \sin \left(\frac{k \pi}{2 K}\right)\right\}.$ to partition the object space into $K$ subregion $\Omega_k, k=1, \ldots K$ and randomly sample the vector $u_k$ of subregion $\Omega_k$. However, how do you take a partitioned random sample in any dimensional $J$ space? A good partition must satisfy:
$$
\cup_{k=1}^J\left(\Omega_k\right)=S^J \text { and } \Omega_{k^{\prime}} \cap_{k^{\prime} \neq k} \Omega_k=\emptyset
$$
With $\mathcal{S}^J=\left\{\lambda \in \mathbb{R}_{>0}^J: \sum_j \lambda_j=1\right\}$ is the set feasible values of random variable $u$. Based on \cite{Das}, we define subregions $\Omega_i$ by points $u=\left(u_1, \ldots, u_K\right) \in U \subset \mathcal{S}^J$ such that:
$$
u_1 \in\{0, \delta, 2 \delta, \ldots, \ldots, 1\} \text { s.t } \frac{1}{\delta}=n \in \mathbb{N}^*
$$

If $1<k<J-1, m_k=\frac{\delta}{u_k}, 1 \leq k<k^\prime-1$, we have:
$$
u_{k^\prime} \in\left\{0, \delta, \ldots,\left(K-\sum_{k=1}^{k^\prime-1} m_k\right) \delta\right\}, u_J=1-\sum_{k=1}^{J-1} u_k
$$

\begin{figure}
    \centering
    \includegraphics{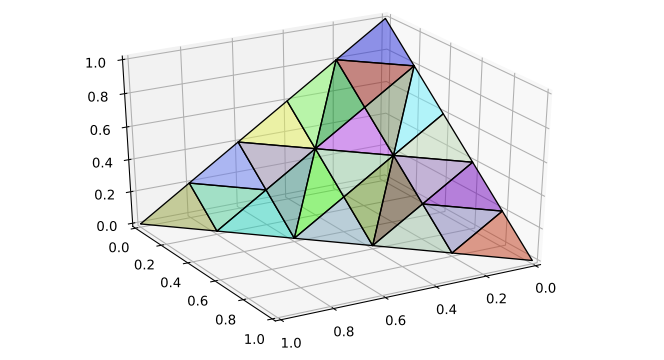}
    \caption{Partitioning algorithm with $J=3,\delta=0.2$}
    \label{partition}
\end{figure}

Using the Delaunay triangulation algorithm \cite{partition} for these points, we obtain the required partition. Figure \ref{partition} is the illustrative result of the algorithm in $3 \mathrm{D}$ space.

\subsection{The update algorithm}
\begin{remark}\label{find_initial}
We propose the following approach to find an acceptable initialization point:
\begin{itemize}
\item \textbf{Step 1:} Randomly initialize a point $x_0 \in \mathbb{R}^n, \eta_r \in (0,1)$.

\item \textbf{Step 2:} If $x_0$ is feasible, stop. Otherwise, proceed to \textbf{Step 3}.

\item \textbf{Step 3:} We can find a descent direction $s_r$ as follows:
$$
\left(s_r, \alpha_r\right)=\arg \min _{s \in \mathbb{R}^n, \alpha \in R} \alpha+\frac{1}{2}\|s\|^2, \text { subject to } \nabla \mathcal{G}_p\left(x_0\right)^T s \leq \alpha
$$

\item \textbf{Step 4:} Update $x_0=x_0+\eta_r s_{r}$ and return to \textbf{Step 2}.
\end{itemize}
\end{remark}

Once an acceptable initialization point is found, we proceed to the proposed algorithm with preference vectors. We can find a descent direction to solve \eqref{P3} applying a similar schema as in problem \ref{projected}.

\begin{equation}\label{3.8}
\begin{cases}\min & g(t, s)=t \\ \text { subject to } & \nabla F_j(x)^T s+\frac{1}{2} s^Ts  -t \leq 0 \quad(1 \leq j \leq m)\\
&  \nabla \mathcal{G}_p(x)^{\top} s +\frac{1}{2} s^Ts  -t \leq 0 \quad (p \in I_\epsilon(x))\\ & (t, s) \in \mathbb{R} \times \mathbb{R}^n.\end{cases}\tag{$P_2$}
\end{equation}
By apply KKT condition to solve problem \ref{3.8}, we have the fomula of descent direction is:
\begin{equation}\label{descent_pref}
   s(x^k)=-\left(\sum_{j=1}^m \lambda^{k}_{j} \nabla F_j(x^{k})+\sum_{p\in I_\epsilon(x^k)} \gamma_p^k\nabla \mathcal{G}_p(x^k)\right),
\end{equation}
where $I_\epsilon(x^k) = \{p|\mathcal{G}_p(x^k)\geq -\epsilon, p = 1,\ldots, K\}$. So we have the following proposed algorithm with preference vectors.
\begin{algorithm}[H]
\caption{Adaptive multi-gradient with preference vectors}
\label{algo_nocons_2}
\textbf{Step 1}. Set $\kappa \in[ 0,1]$, $\sigma \in [0,1]$, $\alpha_1 \in (0, 1)$, and $x^{1}_r \in \mathbb{R}^n$. Initialize $k:=1$.

\textbf{Step 2.} Obtain an acceptable starting point $x^1$ from $x^1_r$.

\textbf{Step 3.} (Main loop)

(a) Find the descent direction by solving the following problem:

\begin{equation} \label{algo2_prob}
    \begin{aligned}
& \min _{\lambda_j, \gamma_p}    ||\sum_{j=1}^m \lambda_j^k\nabla F_j(x^k) + \sum_{p\in I_\epsilon(x^k)} \gamma_p^k\nabla \mathcal{G}_p(x^k)||^2 \\
&\text {subject to } \lambda_j \geq 0, \gamma_p \geq 0, \,\,\sum_{j=1}^m \lambda_j^k+\sum_{p\in I_\epsilon(x^k)}\gamma_p^k=1
\end{aligned}
\end{equation}

With,  $I_\epsilon(x^k):=\{p\in {1,\ldots,K}|\mathcal{G}_p(x^k)\geq -\epsilon\}$. 

So we can obtain the optimal solutions:

$\mathbf{\lambda^k} = (\lambda_j^k)$, $\mathbf{\gamma^k} = (\gamma_p^k)$, $\,\,s(x^k)=-\left(\sum_{j=1}^m \lambda^{k}_{j} \nabla F_j(x^{k})+\sum_{p \in I_\epsilon(x)} \gamma_p^k\nabla \mathcal{G}_p(x^k)\right)$ and $\Theta(x^{k}) =-\dfrac{1}{2}||s(x^k)||^2.$ 

(b) If $\Theta(x^{k})=0$ then \textbf{End}. Otherwise, continue to \textbf{Step 3(c)}

(c) Set $x^{k+1}:=x^{k}+\alpha_{k} s(x^{k}).$

(d) Compute the step size:  

If \begin{equation}\label{kiemtra}
F(x^{k+1}) \leq F(x^{k}) +\sigma \left\langle JF(x^k),x^{k+1}-x^k \right\rangle
\end{equation}

then $\alpha_{k+1} = \alpha_k$, otherwise set $\alpha_{k+1}:=\kappa \alpha_k$.

\textbf{Step 4.} Set $k:=k+1$ and go to \textbf{Step 3.}
\end{algorithm}

\subsection{Convergence of adaptive multi-gradient with preference vector}\label{scalarize}

In this subsection, we show that the Algorithm \ref{algo_nocons_2} can be reformulated as a linear scalarization of tasks with adaptive weight assignment. We first tackle the unconstrained case. Suppose we do not decompose the multi-objective problem and hence remove all constraints from the problem \ref{P3}, it will immediately reduce to the update rule in Lemma \ref{lemma4.2}. It is straightforward to rewrite the corresponding \ref{MOP} with $U=\mathbb{R}^n$ into a linear scalarization form:
\begin{equation}\label{scalarization}
    F(x^k)=\sum_{j=1}^m \lambda_j^k F_j(x^k)
\end{equation}
where we adaptively assign the weights $\lambda_i$ by solving the following problem in each iteration:
$$
\min _{\lambda_j}\frac{1}{2}\left\|\sum_{j=1}^m \lambda_j \nabla F_j\left(x^k\right)\right\|^2, \quad \text { s.t. } \sum_{j=1}^m \lambda_j=1, \quad \lambda_j \geq 0, \forall j=1, \ldots, m .
$$

In the constrained case, we have extra constraint terms $\mathcal{G}_p\left(x^k\right)$. If $\mathcal{G}_p\left(x^k\right)$ is inactivated, we can ignore it. For an activated $\mathcal{G}_p\left(x^k\right)$, assuming the corresponding reference vector is $\boldsymbol{u}_i$, we have:
$$
\nabla \mathcal{G}_p\left(x^k\right)=\left(\boldsymbol{u}_p-\boldsymbol{u}_i\right)^T \nabla F\left(x^k\right)=\sum_{j=1}^m\left(\boldsymbol{u}_{p j}-\boldsymbol{u}_{i j}\right) \nabla F_j\left(x^k\right) .
$$

Since the gradient direction $s(x^k)$ can be written as a linear combination of all $\nabla F_j\left(x^k\right)$ and $\nabla \mathcal{G}_p\left(x^k\right)$ as in \eqref{descent_pref}, the general Pareto MTL algorithm can be rewritten as:
\begin{equation}\label{new_scala}
    F(x^k)=\sum_{j=1}^m \alpha_j F_j\left(x^k\right), \text { where } \alpha_j=\lambda_j+\sum_{p \in I_{\epsilon}(x^k)} \gamma_p\left(\boldsymbol{u}_{pj}-\boldsymbol{u}_{ij}\right),
\end{equation}
where $\lambda_j$ and $\gamma_p$ are obtained by solving problem \eqref{algo2_prob} with assigned reference vector $\boldsymbol{u}_i$.
So basically, we can view the problem of solving \eqref{P3} as simultaneously solving multiple unconstrained multi-objective optimization subproblems. The scalarization form \eqref{scalarization} of the multi-objective problem \ref{MOP} is transformed into \eqref{new_scala} with the same functions $F_j$, differing only in the scalarization coefficients $\alpha_j$. So similar to the proof above, we have some following theorems

\begin{theorem}\label{theorem4_nocons}
		Assume that the sequence $ \left\{ x^k\right\} $ is generated by Algorithm \ref{algo_nocons_2}. Then, each limit point (if any) of the sequence $ \left\{ x^k\right\} $ is a Pareto stationary point of the problem \eqref{MOP} where $U=\mathbb{R}^n$. Define
$
T=\left\{x \in \mathbb{R}^n: F(x) \leq F\left(x^k\right), \forall k\right\} .
$
		\begin{itemize}
			\item[$ \bullet $] if $ F $ is quasiconvex on $\mathbb{R}^n$ and $T \neq \emptyset$, then the sequence $ \left\{x^k \right\} $ converges to a Pareto stationary point of the problem.
			\item[$ \bullet $]  if $ F $ is pseudoconvex on $\mathbb{R}^n$ and $T \neq \emptyset$, then the sequence $ \left\{x^k \right\} $ converges to a weakly efficient solution of the problem.
		\end{itemize}
	\end{theorem}
\section{Numerical experiments}

The configuration of hyperparameters will be presented in each example. In this paper, we will experiment with the proposed algorithm using baselines and metrics as follows:
\begin{itemize}
    \item \textbf{Baselines:} We compare Algorithm \ref{algo_nocons_2} with the PMTL algorithm \cite{PMTL}.
    \item \textbf{Evaluation metric:} The area dominated by Pareto front
is known as Hypervolume \cite{HyperVolume}. The
higher Hypervolume, the better Pareto front.
\end{itemize}
To fairly compare the experimental results of the algorithms, we use simple preference vectors $[\cos(\frac{k\pi}{2K}), \sin(\frac{k\pi}{2K})]$ with $K$ is number of initial points  to evenly divide the 2D space and employ the Delaunay triangulation algorithm to generate preference vectors for the case of 3D space.
\begin{example} \label{ex1}
Consider a convex multi-objective optimization problem without constraints:
$$
\begin{aligned}
&\text{Min } F(x)=\left\{\frac{1}{25} x_1^2+\frac{1}{100}\left(x_2-\frac{9}{2}\right)^2, \frac{1}{25} x_2^2+\frac{1}{100}\left(x_1-\frac{9}{2}\right)^2\right\} \\
&\text { subject to } x \in \mathbb{R}^n
\end{aligned}
$$
\end{example}
The results obtained by using Algorithm \ref{algo_nocons_1} and Algorithm \ref{algo_nocons_2} to solve Example \ref{ex1} with 10 initialization points are as follows:
\begin{figure}[H]
   \begin{minipage}{0.48\textwidth}
     \centering
     \includegraphics[scale=0.4]{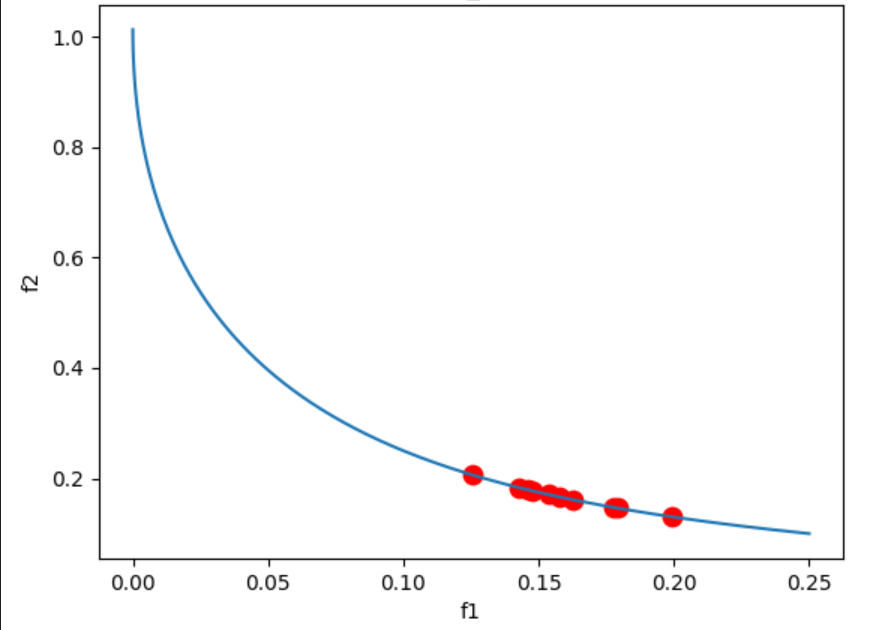}
     \caption{Result of Algo \ref{algo_nocons_1} }\label{Fig:Data1}
   \end{minipage}\hfill
   \begin{minipage}{0.48\textwidth}
     \centering
     \includegraphics[scale=0.4]{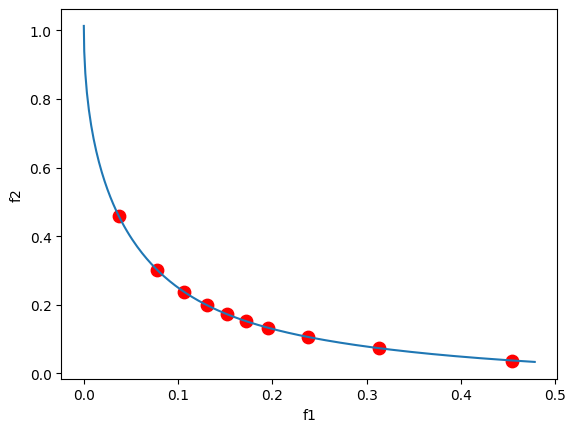}
     \caption{Result of Algorithm \ref{algo_nocons_2}}\label{Fig:Data2}
   \end{minipage}
\end{figure}
For Algorithm \ref{algo_nocons_1}, we can see that the obtained Pareto solutions are concentrated in a region on the Pareto Front. This drawback has been overcome by using constraints to divide the image space in Algorithm \ref{algo_nocons_2}. 

To compare the effectiveness of the proposed algorithm with the PMTL algorithm, we conducted experiments by randomly initializing points to find solutions for Example \ref{ex1}. We have the following result:

\begin{table}[H]
\begin{tabular}{|c|c|c|c|c|c|c|}
\hline
\begin{tabular}[c]{@{}l@{}}Initial \\Point\end{tabular} & \begin{tabular}[c]{@{}l@{}}Iteration \\ Number\end{tabular} & \begin{tabular}[c]{@{}l@{}}Test case \\ number\end{tabular} &\begin{tabular}[c]{@{}l@{}}AVG Time \\ PMTL\end{tabular}  & \begin{tabular}[c]{@{}l@{}}AVG Time \\ Proposed Algo\end{tabular} & \begin{tabular}[c]{@{}l@{}}AVG HV \\ PMTL\end{tabular} & \begin{tabular}[c]{@{}l@{}}AVG Proposed\\ Algo\end{tabular} \\ \hline
50 &450 &8
&1064.64
       & \textbf{915.49}                                                            & 0.99                                                  & \textbf{1.01}       \\

\hline
40 &1000 &7
&629.65
       & \textbf{474.18}                                                            & 1.02                                                  & \textbf{1.03}       \\
\hline
40 &500 &10
&250.18
       & \textbf{230.21}                                                            & 1.01                                                  & \textbf{1.02}       \\
\hline
\end{tabular}
\end{table}

We can see that the proposed algorithm has a faster running time compared to the PMTL algorithm and the obtained solutions converge better to the Pareto Front compared to the PMTL algorithm.
To compare the superiority of the optimization time evenly spread on the Pareto Front, we ran Algorithm \ref{algo_nocons_2} with fewer iterations than the PMTL algorithm and take average results after 10 times.. The results are as follows:

\begin{table}[H]
    \centering
    \renewcommand{\arraystretch}{1.5}
    \begin{tabular}{|c|c|c|c|c|c|c|}\hline
        \multirow{2}{*}{\begin{tabular}[c]{@{}l@{}}Number of \\ Initial Points \end{tabular}} & \multicolumn{3}{c|}{PMTL Algorithm}& \multicolumn{3}{c|}{Algorithm \ref{algo_nocons_2}}\\\cline{2-7}
        & Average HV & Time (s) & Iterations   & Average HV & Time (s) & Iterations \\\hline
        $10$ & $0.97$ & $3.31$ & 200 & $\mathbf{0.97}$ & $\textbf{2.59}$& \textbf{150} \\
        $20$ & $0.96$ & $15.73$ & 200 & $\mathbf{0.96}$ & $\textbf{12.97}$& \textbf{150}\\
        $30$ & $0.95$ & $77.04$ & 300 & $\textbf{0.96}$ & $\textbf{55.07}$ & \textbf{200}\\
        $40$ & $0.95$ & $251.4$ & 300 & $\mathbf{0.95}$ & $\textbf{174.13}$ & \textbf{200}\\
        $50$ & $0.94$ & $835.25$ & 300 & $\textbf{0.95}$ & $\textbf{564.18}$ & \textbf{200}\\\hline
    \end{tabular}
    \caption{Comparison results between Algorithm \ref{algo_nocons_2} and the PMTL algorithm for solving Example \ref{ex1} in various scenarios.}
    \label{table1}
\end{table}
The obtained results show that with fewer iterations, the proposed Algorithm \ref{algo_nocons_2} still achieves an average HV no worse than the PMTL algorithm with more iterations. Consequently, the running speed of Algorithm \ref{algo_nocons_2} is faster than the PMTL algorithm.

\begin{example}\label{ex2}
We consider the following non-convex constrained multi-objective optimization problem:
		$$\begin{array}{ll}\text {Min} & F(x)=\{\dfrac{2x_{1}^{2}+x_{2}^{2}+3}{1+2 x_{1}+8 x_{2}},  \dfrac{x_{1}^{2}+2 x_{2}^{2}+3}{1+8 x_{1} + 2 x_2}\}\\ \text {subject to } & x\in \mathbb{R}^2 \end{array}$$ 
\end{example}
We perform the experiment with Algorithm \ref{algo_nocons_1} and Algorithm \ref{algo_nocons_2} using 10 random initial points, and we obtain the following results:
 \begin{figure}[H]
   \begin{minipage}{0.48\textwidth}
     \centering
     \includegraphics[scale=0.5]{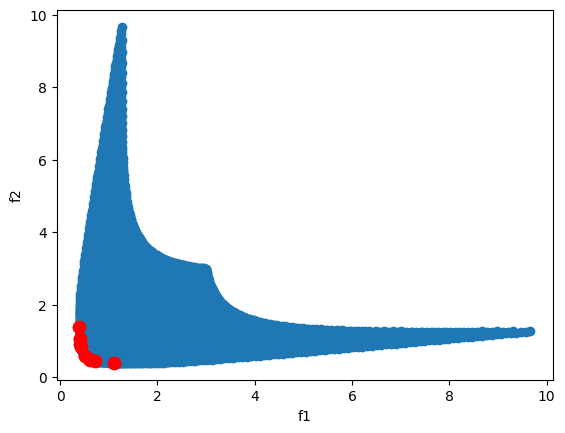}
     \caption{Result of solving Problem 4 of Algorithm \ref{algo_nocons_1}}\label{Fig:Data1}
   \end{minipage}\hfill
   \begin{minipage}{0.48\textwidth}
     \centering
     \includegraphics[scale=0.5]{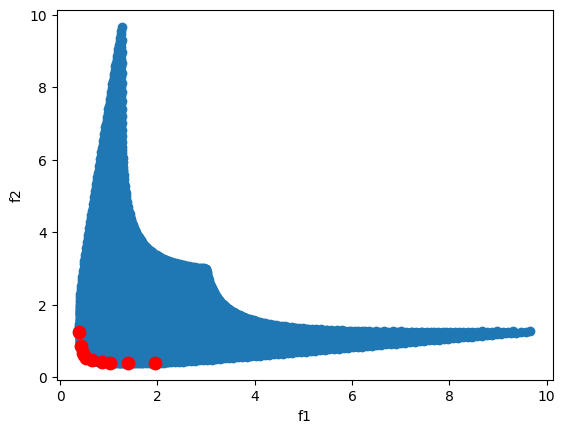}
     \caption{Result of solving Problem 4 of Algorithm \ref{algo_nocons_2}}\label{Fig:Data2}
   \end{minipage}
\end{figure}
We can see that both algorithms are able to find solutions on the Pareto Front, but with the linear constraints that help divide the image space, Algorithm \ref{algo_nocons_2} finds solutions that are more evenly spread out on the Pareto Front. To compare the effectiveness of the proposed algorithm with the PMTL algorithm, we conducted experiments by randomly initializing points to find solutions for Example \ref{ex2}. We have the following result:

\begin{table}[H]
\begin{tabular}{|c|c|c|c|c|c|c|}
\hline
\begin{tabular}[c]{@{}l@{}}Initial \\Point\end{tabular} & \begin{tabular}[c]{@{}l@{}}Iteration \\ Number\end{tabular} & \begin{tabular}[c]{@{}l@{}}Test case \\ number\end{tabular} &\begin{tabular}[c]{@{}l@{}}AVG Time \\ PMTL\end{tabular}  & \begin{tabular}[c]{@{}l@{}}AVG Time \\ Proposed Algo\end{tabular} & \begin{tabular}[c]{@{}l@{}}AVG HV \\ PMTL\end{tabular} & \begin{tabular}[c]{@{}l@{}}AVG Proposed\\ Algo\end{tabular} \\ \hline
35 &1000 &10
&110.04
       & \textbf{94.53}                                                            & 1.14                                                  & \textbf{1.14}       \\

\hline
40 &1000 &5
&152.53
       & \textbf{131.01}                                                            & 1.14                                                  & \textbf{1.14}       \\
\hline

\end{tabular}
\end{table}

We can see that the proposed algorithm has a faster running time compared to the PMTL algorithm and the obtained solutions converge better to the Pareto Front compared to the PMTL algorithm.
To compare the superiority of the optimization time evenly spread on the Pareto Front, we ran Algorithm \ref{algo_nocons_2} with fewer iterations than the PMTL algorithm and take average results after 10 times. The results are as follows:

\begin{table}[H]
    \centering
    \renewcommand{\arraystretch}{1.5}
    \begin{tabular}{|c|c|c|c|c|c|c|}\hline
        \multirow{2}{*}{\begin{tabular}[c]{@{}l@{}}Number of \\ Initial Points \end{tabular}} & \multicolumn{3}{c|}{PMTL Algorithm}& \multicolumn{3}{c|}{Algorithm \ref{algo_nocons_2}}\\\cline{2-7}
        & Average HV & Time (s) & Iterations   & Average HV & Time (s) & Iterations \\\hline
        $10$ & $1.12$ & $9.65$ & 500 & $\mathbf{1.12}$ & $\textbf{5.99}$& \textbf{400} \\
        $20$ & $1.13$ & $21.28$ & 500 & $\mathbf{1.13}$ & $\textbf{11.87}$& \textbf{300}\\
        $30$ & $1.14$ & $50.97$ & 500 & $\textbf{1.14}$ & $\textbf{48.74}$ & \textbf{300}\\
        \hline
    \end{tabular}
    \caption{Comparison results between Algorithm \ref{algo_nocons_2} and the PMTL algorithm for solving Example \ref{ex2} in various scenarios.}
    \label{table1}
\end{table}
The obtained results show that with fewer iterations, the proposed Algorithm \ref{algo_nocons_2} still achieves an average HV no worse than the PMTL algorithm with more iterations. Consequently, the running speed of Algorithm \ref{algo_nocons_2} is faster than the PMTL algorithm.
\begin{example} (Toy example (Lin et al \cite{lin2019pareto})) \label{ex3}
Consider the non-convex unconstrained multi-objective optimization problem:
$$
\begin{aligned}
&\text{Min } F(x)=\left\{1-\exp ^{-\sum_{i=1}^d\left(x_i-\frac{1}{d}\right)^2}, 1-\exp ^{-\sum_{i=1}^d\left(x_i +\frac{1}{d}\right)^2}\right\} \\
&\text { subject to } x \in \mathbb{R}^2
\end{aligned}
$$
\end{example}
We conducted experiments with $d=20$.
We apply Algorithm \ref{algo_nocons_2} to solve the above problem with 10 initial data points, and we obtain the following results:
\begin{figure}[H]
   \begin{minipage}{0.48\textwidth}
     \centering
     \includegraphics[scale=0.4]{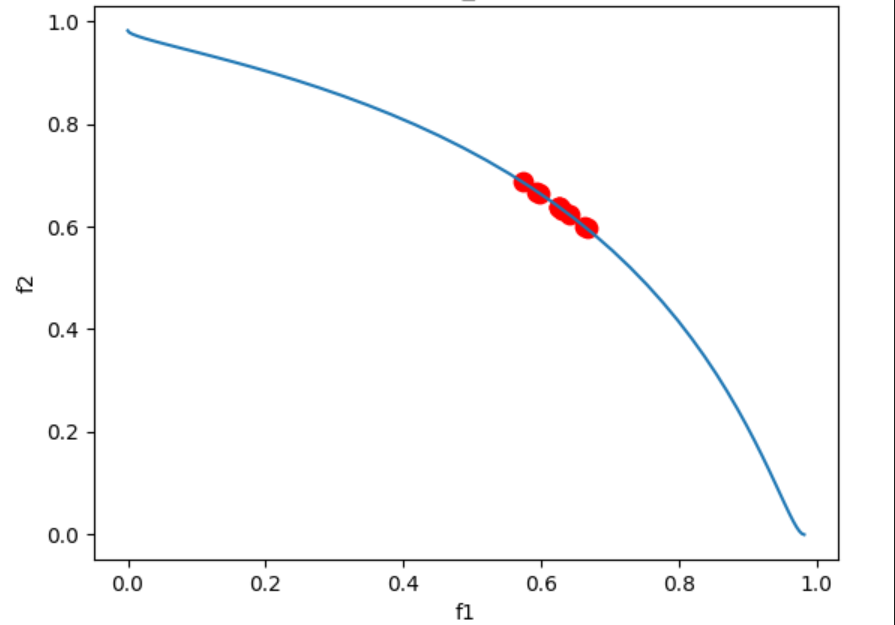}
     \caption{Result of Algo \ref{algo_nocons_1}}\label{Fig:Data1}
   \end{minipage}\hfill
   \begin{minipage}{0.48\textwidth}
     \centering
     \includegraphics[scale=0.4]{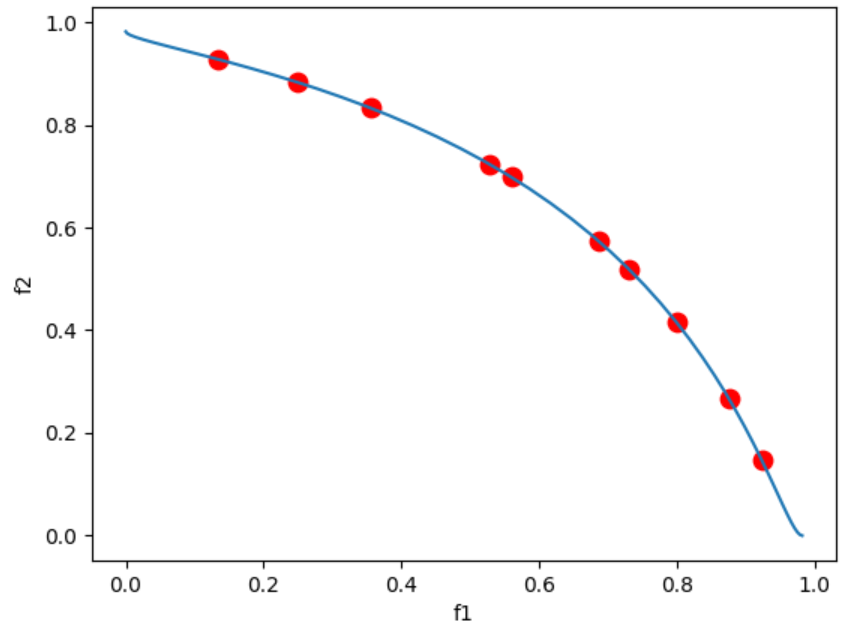}
 \caption{Result of Algo \ref{algo_nocons_2}}
   \end{minipage}
\end{figure}
The solutions obtained by Algorithm \ref{algo_nocons_1} concentrate in a region on the Pareto Front, while using Algorithm \ref{algo_nocons_2}, we obtain solutions that are all on the Pareto Front and spread evenly across the surface. 

Firstly, we compare our proposed Algorithm \ref{algo_nocons_2} with MGDA \cite{desideri} algorithm

\begin{table}[H]
\centering
\begin{tabular}{|c|c|c|c|c|}
\hline
\begin{tabular}[c]{@{}l@{}}Initial \\Point\end{tabular} & \begin{tabular}[c]{@{}l@{}}Iteration \\ Number\end{tabular} & \begin{tabular}[c]{@{}l@{}}Test case \\ number\end{tabular}  & \begin{tabular}[c]{@{}l@{}}AVG HV \\ PMTL\end{tabular} & \begin{tabular}[c]{@{}l@{}}AVG Proposed\\ Algo\end{tabular} \\ \hline
50 &100 &10                                                           & 0.85                                                   & \textbf{1.52}        \\
40 &200 &13                                                           & 0.87                                                   & \textbf{1.49}        \\
45 &150 &15                                                           & 0.88                                                   & \textbf{1.53}        \\
\hline

\end{tabular}
\end{table}

The HV achieved by the proposed algorithm is significantly higher. Utilizing the proposed algorithm will help distribute the obtained solutions evenly across the Pareto front and yield higher HV. Secondly, to compare the effectiveness of the proposed algorithm with the PMTL algorithm, we conducted experiments by randomly initializing points to find solutions for Example \ref{ex3}. We measured the solving time of both algorithms and the Hyper Volume value of the Pareto Front obtained by each algorithm. The results obtained are as follows

\begin{table}[H]
\begin{tabular}{|c|c|c|c|c|c|c|}
\hline
\begin{tabular}[c]{@{}l@{}}Initial \\Point\end{tabular} & \begin{tabular}[c]{@{}l@{}}Iteration \\ Number\end{tabular} & \begin{tabular}[c]{@{}l@{}}Test case \\ number\end{tabular} &\begin{tabular}[c]{@{}l@{}}AVG Time \\ PMTL\end{tabular}  & \begin{tabular}[c]{@{}l@{}}AVG Time \\ Proposed Algo\end{tabular} & \begin{tabular}[c]{@{}l@{}}AVG HV \\ PMTL\end{tabular} & \begin{tabular}[c]{@{}l@{}}AVG Proposed\\ Algo\end{tabular} \\ \hline
50 &200 &40 &264.74        & \textbf{224.06}                                                            & 1.48                                                   & \textbf{1.53}        \\
\hline
50 &300 &10 &245.86        & \textbf{181.84}                                                         & 1.50                                                  &\textbf{ 1.53}        \\
\hline
60 &200 &10 &259.52
        & \textbf{236.51}                                                         & 1.47                                                  &\textbf{ 1.51}        \\
\hline
\end{tabular}
\end{table}

We can see that the proposed algorithm has a faster running time than the PMTL algorithm  and at the same time, the obtained solutions converge better to the Pareto Front than the PMTL algorithm. To compare the superiority of solving evenly-distributed Pareto Fronts with respect to optimization time, we executed Algorithm \ref{algo_nocons_2} with fewer iterations than PMTL. The results demonstrate that with fewer iterations, Algorithm \ref{algo_nocons_2} still achieves an average HV no worse than PMTL with more iterations. Specifically:

\begin{table}[H]
    \centering
    \renewcommand{\arraystretch}{1.5}
    \scalebox{0.9}{ \begin{tabular}{|c|c|c|c|c|c|c|}\hline
        \multirow{2}{*}{\begin{tabular}[c]{@{}l@{}}Number of \\ Initial Points \end{tabular}} & \multicolumn{3}{c|}{PMTL Algorithm}& \multicolumn{3}{c|}{Algorithm \ref{algo_nocons_2}}\\\cline{2-7}
        & Average HV & Time (s) & Iterations   & Average HV & Time (s) & Iterations \\\hline
        $10$ & $1.49$ & $2.01$ & 100 & $\mathbf{1.49}$ & \textbf{1.37}& \textbf{50} \\
        $20$ & $1.53$ & $4.38$ & 100 & $\mathbf{1.53}$ & $\mathbf{2.45}$& \textbf{80}\\
        $30$ & $1.49$ & $17.97$ & 100 & $\mathbf{1.51}$ & $\mathbf{10.79}$ & \textbf{80}\\
        $40$ & $1.49$ & $54.01$ & 100 & $\mathbf{1.49}$ & $\mathbf{35.73}$ & \textbf{75}\\
        $50$ & $1.44$ & $99.42$ & 100 & $\mathbf{1.47}$ & $\mathbf{67.67}$ & \textbf{75}\\\hline
    \end{tabular} }
    \caption{Comparison results between Algorithm \ref{algo_nocons_2} and the PMTL algorithm in various scenarios.}
    \label{table1}
\end{table}

\begin{example}
    \label{vd3_no} (Toy example (Lin et al \cite{lin2019pareto}))
Consider a three-dimensional multi-objective optimization example with the following non-convex functions:
$$
\begin{aligned}
& \text{Min } F(x) =&\left\{1-\exp ^{-\sum_{i=1}^d\left(x_i-\frac{1}{d}\right)^2}, 1-\exp ^{-\sum_{i=1}^d\left(x_i +\frac{1}{d}\right)^2}, \right.\\
& &\left.
 2-\exp ^{-\sum_{i=1}^d\left(x_i-\frac{1}{d}\right)^2} -\exp ^{-\sum_{i=1}^d\left(x_i +\frac{1}{d}\right)^2}\right\} \\
&\text { s.t } x \in \mathbb{R}^d &
\end{aligned}
$$
\end{example} \

We proceed to solve Example \ref{vd3_no} with dimensionality $d=20$. Applying Algorithm \ref{algo_nocons_2} to solve the example with 10 initial points, we obtain the following results:
\begin{figure}[H]
\centering
\includegraphics{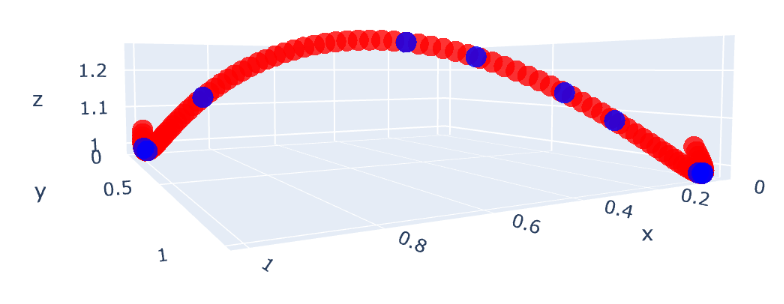}
\caption{Results of Algorithm \ref{algo_nocons_2} solving Example \ref{vd3_no}}
\label{fig:enter-label}
\end{figure}

The Pareto solutions found by Algorithm \ref{algo_nocons_2} are evenly distributed on the Pareto Front.

In summary, the solutions found by the proposed algorithm converge better to the Pareto Front than PMTL. The proposed algorithm can be seen as an improved version of PMTL with the ability to search for Pareto solutions more evenly and faster on unconstrained multi-objective optimization problems.
\section{Applications to multi-task learning}
\subsection{Proposed Algorithm}
Similarly to \cite{lin2019pareto}, we will use $K$ preference vectors ${\mathbf{u_1}, \mathbf{u_2},\ldots,\mathbf{u_K}}$ to represent the trade-off preferences among tasks in the MTL problem. Additionally with $\theta$ are model's parameters, we incorporate linear inequality constraints $\mathcal{G}_p\left(\theta_t\right)=\langle\boldsymbol{u}_p-\boldsymbol{u}_k,\mathcal{L}\left(\theta_t\right)\rangle \leq 0, \forall p=1, \ldots, K$ to partition the image space into subproblems. Then, we propose using an adaptive step size to update the learning rate of the Neural network. We need to solve the following problem
\begin{equation}
    \begin{aligned}
& \min _\theta \mathcal{L}(\theta)=\left(\mathcal{L}_1(\theta), \mathcal{L}_2(\theta), \cdots, \mathcal{L}_m(\theta)\right)^{\mathrm{T}} \\
& \text { s.t. } \quad \mathcal{G}_j\left(\theta_t\right)=\left(\boldsymbol{u}_j-\boldsymbol{u}_k\right)^T \mathcal{L}\left(\theta_t\right) \leq 0, \forall j=1, \ldots, K,
\end{aligned}\tag{$P_3$}
\end{equation}

By applying KKT condition and similar to transformation in Section \ref{scalarize}, we get:
 $$\begin{aligned}
     s(\theta_t) &=-\left(\sum_{j=1}^m \lambda^{k}_{j} \nabla \mathcal{L}_j(\theta_t)+\sum_{p \in I(\theta_t)} \gamma_p^k\nabla \mathcal{G}_p(\theta_t)\right) \\
     \mathcal{L}(\theta_t)&=\sum_{j=1}^m \alpha_j \mathcal{L}_j\left(\theta_t\right), \text { where } \alpha_j=\lambda_j+\sum_{p \in I_{\epsilon}(\theta_t)} \gamma_p\left(\boldsymbol{u}_{pj}-\boldsymbol{u}_{ij}\right)
 \end{aligned}$$
where $I_\epsilon(\theta_t) = \{p|\mathcal{G}_p(x^k)\geq -\epsilon, p = 1,\ldots, K\}$. The detailed algorithm for MTL is presented below.
\begin{algorithm}[H]
\caption{Proposed Algorithm for Solving MTL Problem}
\label{algo_PMTL}
\textbf{Input:} Set of priority vectors $\{\mathbf{u_1}, \mathbf{u_2},\ldots,\mathbf{u_K}\}$.
\begin{algorithmic}[1]
\For{$k = 1$ to $K$}
    \State Set $\kappa \in [0,1]$, $\sigma \in [0,1]$, $\alpha_1 \in (0, +\infty)$.
    \State Initialize parameter set $\theta_r^{k}$ for the Neural network.
    \State Find effective parameters $\theta_0^{k}$ from $\theta_r^{k}$ using the descent method.
    \For{$t =1$ to $T$}
            \State Find the descent direction by solving the following problem:
            
            $$
            \begin{aligned}
            & \min _{\lambda_j, \gamma_p}    ||\sum_{j=1}^m \lambda_j^k\nabla \mathcal{L}_j(\theta^k_t)^T + \sum_{p \in I(\theta^k_t)} \gamma_p^k\nabla \mathcal{G}_p(\theta^k_t)||^2 \\
            &\text { s.t. } \lambda_j \geq 0, \gamma_p \geq 0, \,\,\sum_{j=1}^m \lambda_j^k+\sum_{p \in I(\theta^k_t)} \gamma_p^k=1
            \end{aligned}
            $$
            
            \State Where $I_\epsilon(\theta^k_t):=\{p\in {1,\ldots,K}|\mathcal{G}_p(\theta^k_t)\geq -\epsilon\}$.
            \State Determine the descent direction:
            
            $$s(\theta^k_t)=-\left(\sum_{j=1}^m \lambda^{k}_{j} \nabla \mathcal{L}_j(\theta^{k}_t)+\sum_{p \in I(\theta^k_t)} \gamma_p^k\nabla \mathcal{G}_p(\theta^k_t)\right)$$
            
            \State Set $\theta^{k}_{t+1}:=\theta^{k}_t+\alpha_{t} s(\theta^{k}_t).$
            \State Compute step size:  
            \State If  $\mathcal{L}(\theta^{k}_{t+1}) \leq \mathcal{L}(\theta^{k}_t) + \sigma \left\langle s(\theta^k_t),\theta^{k}_{t+1}-\theta^k_t \right\rangle$ 
            \State then $\alpha_{t+1} = \alpha_t$, otherwise set $\alpha_{t+1}:=\kappa \alpha_t$.
            \State Update iteration step $t:=t+1$.
    \EndFor
\EndFor
\end{algorithmic}
\textbf{Output: }  The set of solutions for all subproblems with different trade-offs $\{\theta^{k}_T|k=1,\ldots,K\}$.
\end{algorithm}

\subsection{Applications to Computer Vision}
In this section, we will implement the proposed algorithm and evaluate its effectiveness on a two-task image classification problem using the Multi-MNIST dataset provided at \cite{lin2019pareto}. At the same time, we will compare the effectiveness of the proposed algorithm with the PMTL algorithm.

Multi-MNIST is a dataset constructed by combining two digits from the MNIST dataset \cite{lecun} into one image. In the MTL problem of image classification, the model simultaneously learns two loss functions to predict two digits, one in the top-left corner and one in the bottom-right corner, nested within each other. The training set consists of 120,000 images, and the test set consists of 20,000 images. An illustrative image of the dataset is shown below:
\begin{figure}[H]
    \centering
    \includegraphics[scale = 0.7]{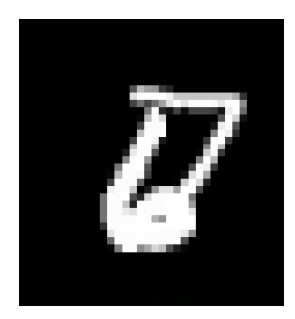}  \label{fig:my_label}
\end{figure}

We construct a LeNet model similar to the one in \cite{lecun}. Our two tasks are image classification tasks, so the objective or loss function in our model is represented as Cross-entropy loss:
$$
\begin{aligned}
    L_1(\hat{y}_1, y_1) &=-\sum_k^K y^{(k)}_1 \log \hat{y}^{(k)}_1 \\
    L_2(\hat{y}_2, y_2) &=-\sum_k^K y^{(k)}_2 \log \hat{y}^{(k)}_2
\end{aligned}
$$
Here, $y^{(k)}_i, (i=1,2)$ takes on values of 0 or 1, representing the correct label $k$ that has been correctly classified, and $\hat{y}^{(k)}_i, (i=1,2)$ are the predicted label values of the model. We need to simultaneously learn both Cross-entropy loss functions. We carry out the training process with 10 reference vectors $[\cos(\frac{k\pi}{2K}), \sin(\frac{k\pi}{2K})],k=1,2,\ldots,10$ for every 100 epochs. We illustrate the results with the reference vector $u_1=(1,0)$ as follows:
\begin{figure}[H]
    \centering
{{\includegraphics[scale=0.2]{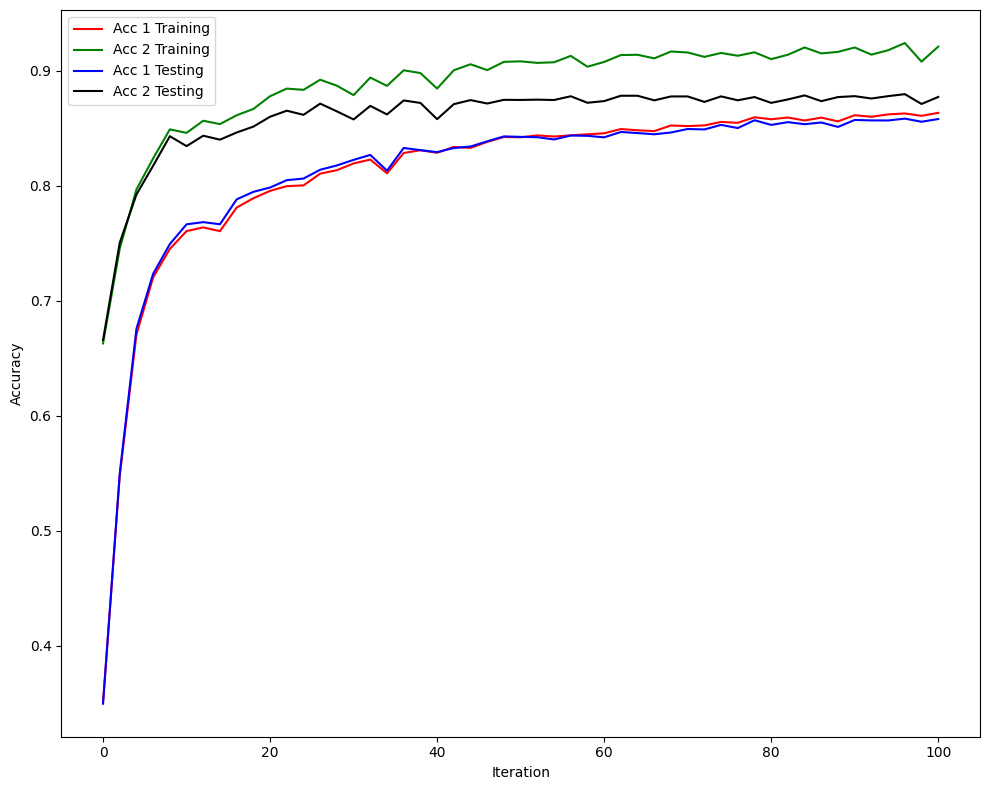} }}%
    \quad
    {{\includegraphics[scale=0.2]{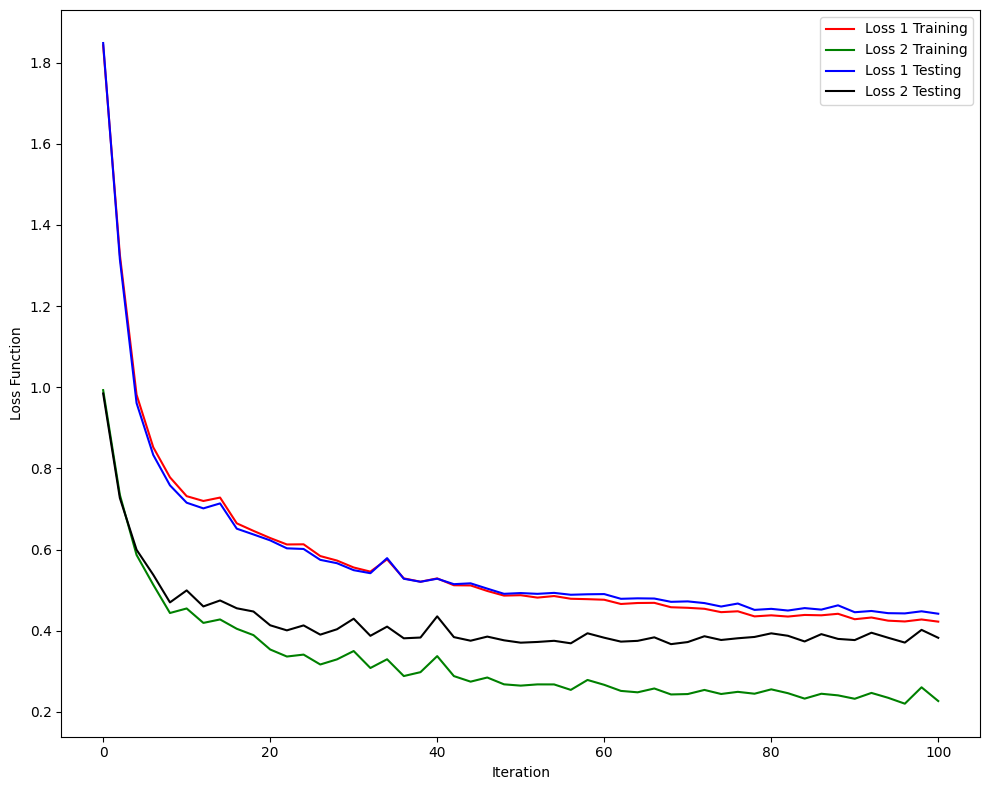} }}
    \caption{The test result}
    \label{fig:example}
\end{figure}
The experimental results show that both image classification tasks are learned simultaneously and relatively effectively. The table comparing the accuracy values on the Test set between the proposed algorithm and PMTL will be illustrated below: 
\begin{figure}[h]
  \begin{minipage}[h]{0.45\textwidth}
    \includegraphics[width=\textwidth]{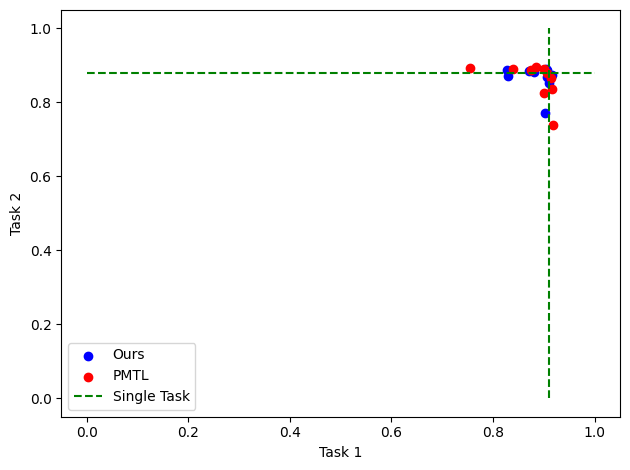}
  \end{minipage}
\hfill
  \begin{minipage}[h]{0.5\textwidth}
    \centering
    \begin{tabular}{|c|c|c|}
      \hline
     &\begin{tabular}[c]{@{}l@{}}Average Accuracy \\Task 1\end{tabular}   &  \begin{tabular}[c]{@{}l@{}}Average Accuracy \\Task 2\end{tabular} \\
      \hline
      PMTL & $0.87 \pm 0.01$ & $0.85 \pm 0.01$ \\
      Ours & \textbf{$0.88 \pm 0.01$} &  \textbf{$0.86 \pm 0.01$}\\
      \hline
    \end{tabular}
  \end{minipage}

\end{figure}

We can see that the effectiveness of the proposed algorithm is quite similar to the PMTL algorithm. The comparison of the training process between the proposed algorithm (blue line for task 1 and black line for task 2) and the PMTL algorithm (red line for task 1 and green line for task 2) for the reference vector $\mathbf{u_1}$ is shown in the following image:
\begin{figure}[H]
    \centering
{{\includegraphics[scale=0.2]{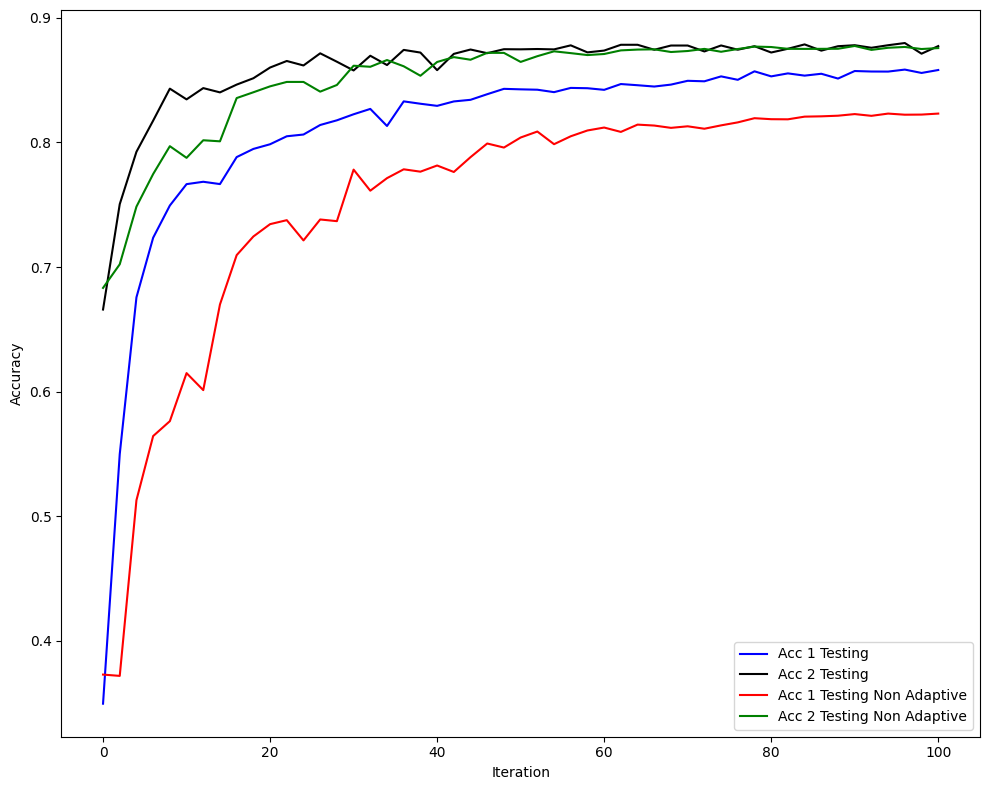} }}%
    \,\,
    {{\includegraphics[scale=0.2]{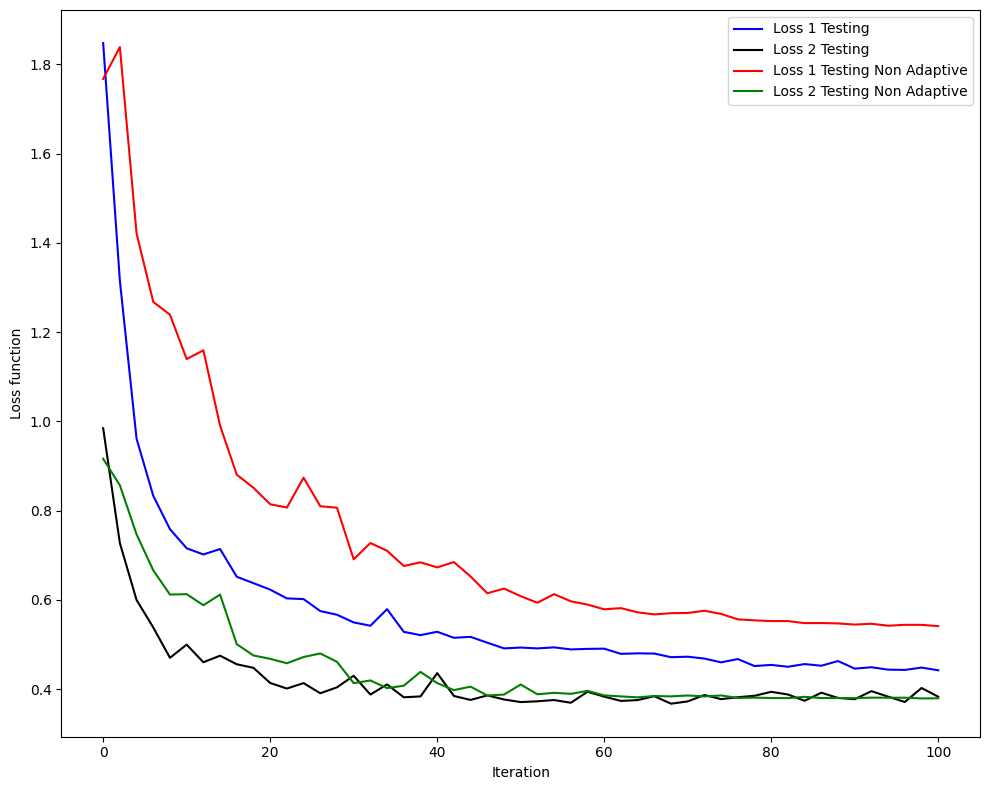} }}
    \caption{Compare the proposed algorithm and PMTL}
    \label{fig:example}
\end{figure}
The proposed algorithm and PMTL algorithm are two multi-task learning methods for image classification tasks. The experimental results show that both algorithms can learn two tasks simultaneously with relatively similar effectiveness, as seen in the comparison of accuracy values on the Test set. However, in terms of the training process, the proposed algorithm has an advantage over PMTL in that it uses an Adaptive step size, which leads to a more consistent decrease in the loss function and therefore a continuous increase in accuracy during training. Meanwhile, PMTL uses a fixed step size, which may lead to slower convergence or oscillation during training. In summary, while both algorithms can achieve similar performance in multi-task learning for image classification, the proposed algorithm has the advantage of a more stable and consistent training process thanks to its Adaptive step size.
\subsection{Application to Natural Language Processing}
Our inquiry focused exclusively on the Drug Review dataset, as outlined in the study conducted by Graßer et al. in 2018. The dataset consists of user assessments of specific prescriptions, information about relevant diseases, and a user rating indicating overall happiness.
We analyze two tasks: (1) predicting the drug's rating using regression and (2) classifying the patient's state.
The dataset comprises 215,063 samples. $90\%$ of the data is utilized, while situations lacking adequate user feedback are eliminated. Next, there are 100 condition labels and 188155 samples. 

\begin{figure}[H]
    \centering
    \includegraphics[scale = 0.7]{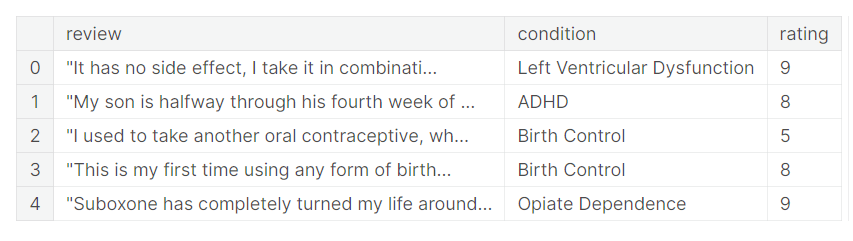}

    \label{fig:enter-label}
\end{figure}
After preprocessing the customer reviews, we proceed to embed the reviews using the GLOVE method \cite{glove}. This conversion transforms the words and sentences into numerical matrix space. Subsequently, we pass this embedded matrix through a TextCNN model \cite{textcnn}. 
We construct the loss function using two Cross-Entropy functions for the multi-label classification problem:
$$
\begin{aligned}
    L_1(\hat{y}_1, y_1) &=-\sum_k^K y^{(k)}_1 \log \hat{y}^{(k)}_1 \\
    L_2(\hat{y}_2, y_2) &=-\sum_k^K y^{(k)}_2 \log \hat{y}^{(k)}_2
\end{aligned}
$$
Here $y^{(k)}_i, (i=1,2)$  takes values of 0 or 1, representing whether label 
$k$ is correctly classified, and $\hat{y}^{(k)}_i, (i=1,2)$ are the predicted label values by the model.

We proceed to train the model using Algorithm \ref{algo_PMTL} with 5 preferrence vectors  $[\cos(\frac{k\pi}{2K}), \sin(\frac{k\pi}{2K})],k=1,2,\ldots,5$ where each priority vector goes through 50 epochs. The hyperparameters, after tuning, are set to 
 $\sigma = 0.05, \kappa = 0.95$ and learning rate $lr = 0.001$. The Train/Test ratio is $3:1$. The results of the training process with loss function metrics and accuracy on the Test set are as follows:
\begin{figure}[H]
    \centering
{{\includegraphics[scale=0.2]{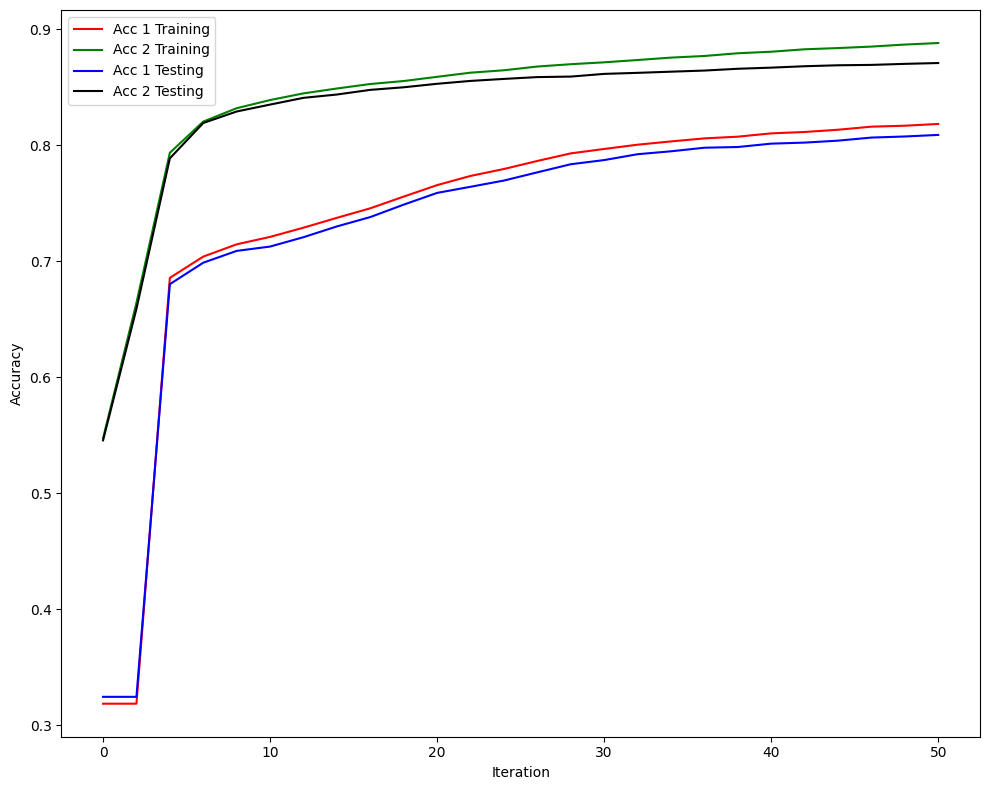} }}%
    \quad
    {{\includegraphics[scale=0.2]{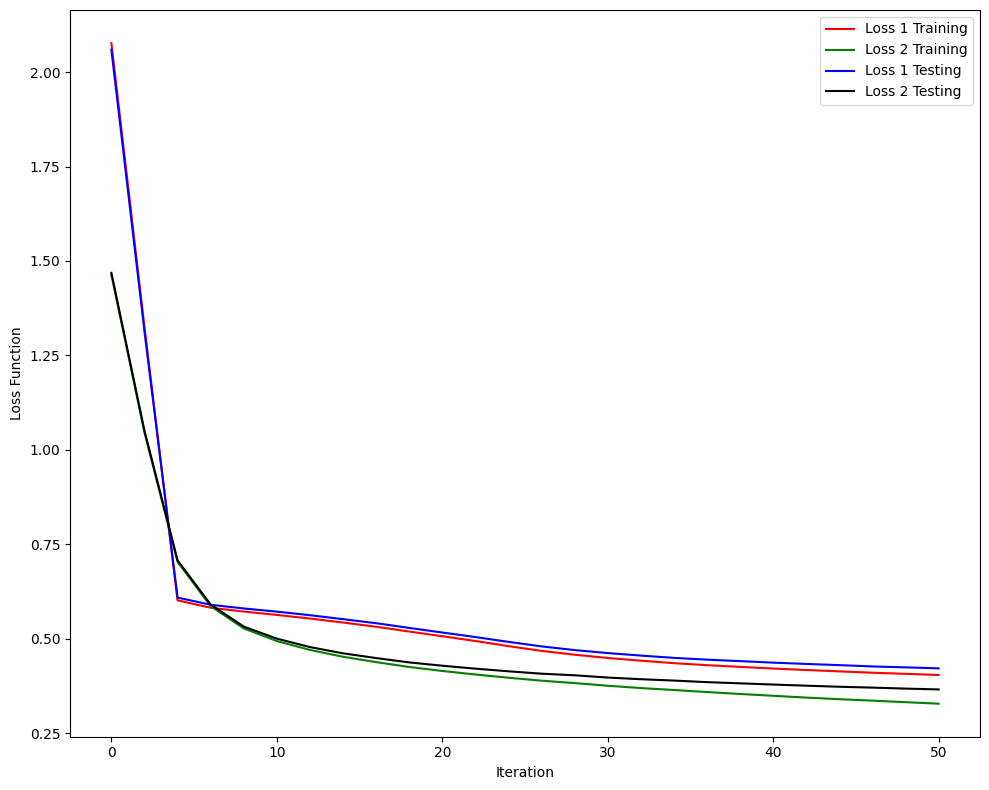} }}
    \caption{Training progress of Algorithm \ref{algo_PMTL} on the Drug Review dataset.}
    \label{fig:example}
\end{figure}
We observe that the loss function values decrease steadily for both the Train and Test sets, and the accuracy consistently increases. This indicates that the model effectively learns both tasks simultaneously.
To compare the performance of the proposed algorithm with the PMTL algorithm, keeping 
$lr = 0.0001$ (the locally optimal hyperparameter after tuning) unchanged for a fair evaluation, we obtain the following results:

\begin{figure}[H]
\centering
{{\includegraphics[scale=0.2]{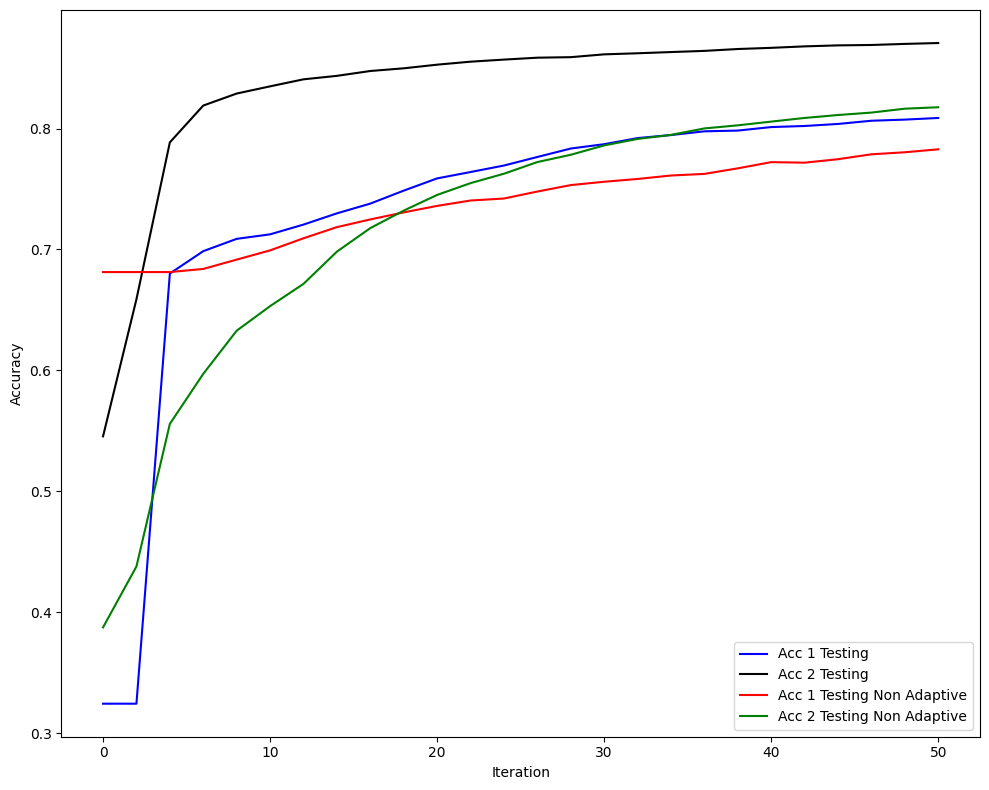} }}%
\quad
{{\includegraphics[scale=0.2]{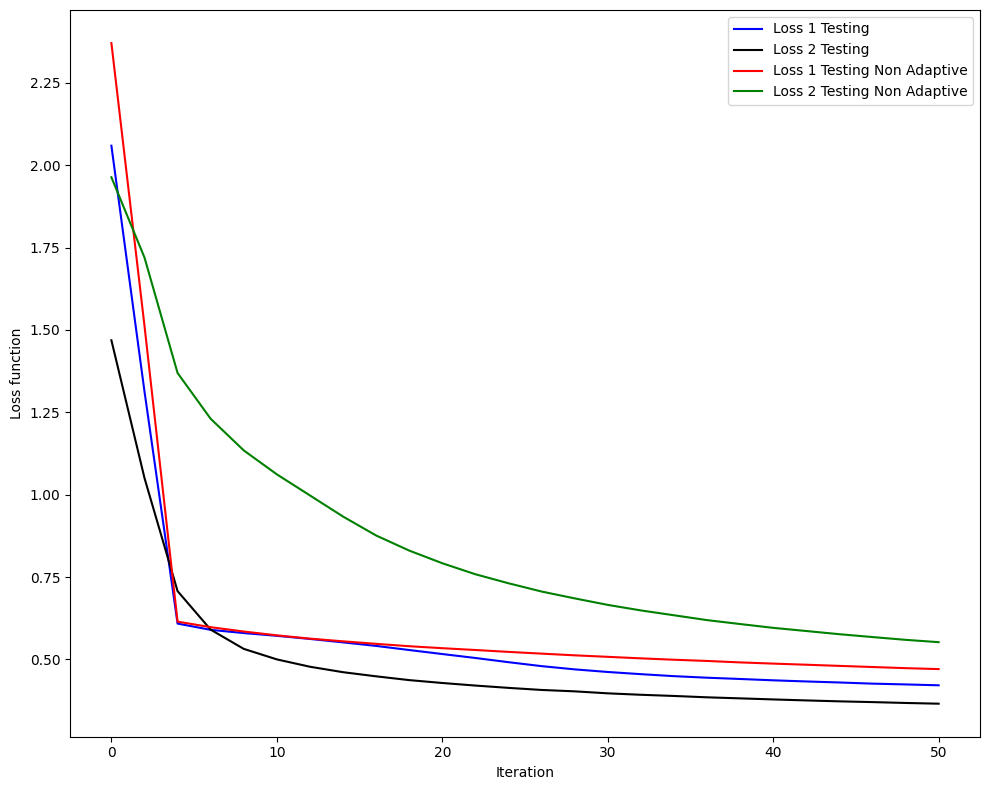} }}
\caption{Comparison of Algorithm \ref{algo_PMTL} and PMTL when run on Drug Review dataset}
\label{fig:example}
\end{figure}

The results show that after only 20 epochs, the accuracy of Algorithm \ref{algo_PMTL} is significantly higher than that of the PMTL algorithm. This is the result of appropriately adjusting the learning rate during training based on the adaptive update condition. We then measure the average accuracy results of the two PMTL algorithms and Algorithm \ref{algo_PMTL} over 10 different runs with the same set of initialization parameters, and the obtained results are as follows:
\begin{minipage}[h]{0.45\textwidth}
  \begin{figure}[H]
      \centering
      \includegraphics[width=\textwidth]{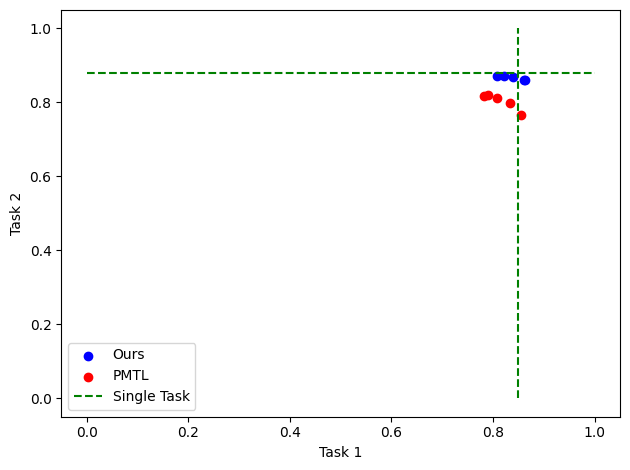}
    \caption{Comparison of the Pareto Fronts between PMTL and Algorithm \ref{algo_PMTL} on the Drug Review dataset.}
  \end{figure}
   
  \end{minipage}
\hfill
  \begin{minipage}[h]{0.5\textwidth}
    \centering
    \begin{table}[H]
        \centering
\begin{tabular}{|c|c|c|}
      \hline
     &\begin{tabular}[c]{@{}l@{}}Average \\  Accuracy Task 1\end{tabular}   &  \begin{tabular}[c]{@{}l@{}}Average \\ Accuracy  Task 2\end{tabular} \\
      \hline \hline
      PMTL & $0.82 \pm 0.01$ & $0.83 \pm 0.01$ \\
      Proposed & $\mathbf{0.84 \pm 0.01}$ &  $\mathbf{0.86 \pm 0.01}$\\
      \hline
    \end{tabular}
        \caption{Comparison between Algorithm \ref{algo_PMTL} and PMTL when run on the Drug Review dataset.} 
        \label{tab:my_label}
    \end{table}

  \end{minipage}

\vspace{ 1 cm}
  The results show the effectiveness of Algorithm \ref{algo_PMTL} compared to PMTL. The Pareto Front found by Algorithm \ref{algo_PMTL} is significantly superior to the Pareto Front found by PMTL.

\section{Conclusion}
We proposed a novel easy adaptive step-size process in a wide family of solution methods for optimization problems with non-convex objective functions. This approach does not need any line-searching or prior knowledge, but rather takes into consideration the iteration sequence's behavior. As a result, as compared to descending line-search approaches, it significantly reduces the implementation cost of each iteration. We demonstrated technique convergence under simple assumptions. We demonstrated that this new process produces a generic foundation for optimization methods. The preliminary results of computer experiments demonstrated the new procedure's efficacy.
	

\end{document}